\documentclass[preprint,12pt]{elsarticle}

\usepackage{amssymb}
\usepackage{amsmath}
\usepackage{amsthm}

\usepackage{tikz}
\usepackage{tikz-3dplot}
\usepackage{pgfplots}
\usepackage{subcaption}
\pgfplotsset{compat=1.6}
\usetikzlibrary{calc,patterns,angles,quotes}
\usepackage{graphicx}
\usepackage{xcolor}
\usepackage{hyperref}

\newcommand{\x}{\boldsymbol{x}}
\newcommand{\y}{\boldsymbol{y}}
\newcommand{\n}{\boldsymbol{n}}

\newcommand{\nx}{\boldsymbol{n}(\boldsymbol{x})}

\newcommand{\diff}{\mathrm{d}}
\newcommand{\V}{\boldsymbol{V}}

\theoremstyle{thmstylethree}%
\newtheorem{definition}{Definition}%

\journal{Engineering Analysis with Boundary Elements}

\begin{document}

\begin{frontmatter}



\title{Evaluating singular and near-singular integrals on $C^2$ smooth surfaces with quadratic geometric approximation and closed form expressions}


\author[1]{Andrew Zheng} 
\ead{andrewf.zheng@mail.utoronto.ca}
\author[1]{Spyros Alexakis}
\ead{alexakis@math.toronto.edu}
\author[1]{Adam R. Stinchcombe\corref{cor1}}
\ead{stinch@math.toronto.edu}

\cortext[cor1]{Corresponding author}

\affiliation[1]{organization={Department of Mathematics, University of Toronto},
            addressline={40 St George Street}, 
            city={Toronto},
            postcode={M5S 2E4}, 
            state={Ontario},
            country={Canada}}

\begin{abstract}
Most Fredholm integral equations involve integrals with weakly singular kernels. Once the domain of integration is discretized into flat triangular elements, these weakly singular kernels become strongly singular or near-singular. 
Common methods to compute these integrals when the kernel is a Green's function include the Duffy transform, polar coordinates with closed analytic formulas, and singularity extraction. However, these methods do not generalize well to the normal derivatives of Green's functions due to the strongly singular behavior of these functions on triangular elements. 
We provide methods to integrate both the Green's function and its normal derivative on smooth surfaces discretized by triangular elements in three dimensions for many commonly encountered differential operators. 
For strongly singular integrals involving normal derivatives of Green's functions, we introduce a more refined approximation of the derivatives of the Green's function on flat triangles. This method uses geometric information of the true surface of integration to approximate the original integral on the true domain using push-forward maps. This is better than simply setting the singular integrals to zero, while being faster than adaptive refinement methods. 
We provide an algorithm for explicit computations on triangles, and present necessary analytic formulas that the algorithm requires in the appendix.
\end{abstract}


\begin{highlights}
\item Expressions for singular and near-singular integrals on triangular elements.
\item Distance preserving algorithm to get integration bounds on flat triangles. 
\end{highlights}


\end{frontmatter}


\section{Introduction}\label{sec: introduction}

Let $\Omega\subset \mathbb{R}^3$ with $C^2$ boundary. Consider the problem of integrating a function $f$ over $\partial\Omega$. Assuming that $f\in L^1(\partial\Omega)$, this can be a difficult problem numerically. 
Unless there are simple analytic formulas for the integral, one will have to create a quadrature scheme for this integral.
As the geometry of $\partial\Omega$ becomes more complicated, an accurate quadrature scheme will be harder to find.

One common algorithm to integrate functions on arbitrary smooth boundaries is to first discretize $\partial\Omega$ into elements on which very accurate quadrature schemes exist \cite{hsiao2008boundary}. In 2D, common elements are line segments, while in 3D they are simplexes or triangles which we denote by $\Delta$. The method of approximating $\partial\Omega$ with a finite number of triangles is often referred to as triangulation, and we denote by $\partial\Omega_{\Delta}$ as the collection of those triangles. 
There are many algorithms available that can perform this triangulation \cite{persson2004simple, delaunay1934spheres} and simple implementations can be found in most commonly used programming languages. 

After approximating $\partial\Omega$ with $\partial\Omega_{\Delta}$, we can write
\begin{equation}
    I(f) = \int_{\partial\Omega}f(\x)\,\diff S(\x) \approx \sum_{\Delta \in \partial\Omega_{\Delta}} \int_{\Delta} f(\x)\, \diff S(\x).
\end{equation}
We use the notation that bold-face variables $\x$ represent vectors, while $x$ represents real values. The function $f$ is often weakly singular on $\partial\Omega$ and strongly singular on the triangles $\Delta$. Thus, by triangulation, we approximate a finite integral with an undefined integral.

A common example of a function that is weakly integrable on $\partial\Omega$ but not on its triangulation is the normal derivative of Green's function for the Laplace equation in three dimensions, which shows up in many electrostatics problems \cite{griffith2013electrophysiology}. The Green's function and its normal derivative in $\x$ are 
\begin{equation}
    G_{\mathbb{R}^3}(\x, \y) = -\frac{1}{4\pi |\x - \y|}, \label{eq: Green in 3D} \quad
    K_{\mathbb{R}^3}(\x, \y) = \frac{\partial G_{\mathbb{R}^3}}{\partial \nx}(\x, \y) = \frac{(\x-\y)\cdot \nx}{4\pi |\x-\y|^3}.
\end{equation}

In more generality, the integrals we wish to evaluate are
\begin{align}
    \int_{\Delta} G(\x, \y)p(\y)\, \diff S(\y) &=  \int_{\Delta} \frac{g(\x, \y)}{|\x-\y|}p(\y)\, \diff S(\y), \label{eq: general integral problem green} \\
    \int_{\Delta} \frac{\partial G}{\partial \nx}(\x, \y)p(\y)\, \diff S(\y) &= \int_{\Delta} \frac{h(\x, \y) (\x - \y)\cdot \nx}{|\x-\y|^3}p(\y)\, \diff S(\y), \label{eq: general integral problem gradient green}
\end{align}
in which $G$ has a type $r^{-1} = |\x-\y|^{-1}$ singularity, $g, h$ are smooth functions, and $\Delta$ is some triangle in $\mathbb{R}^3$. 
Drawing inspiration from potential theory \cite{atkinson1997numerical}, $\x$ is often called the target and $\Delta$ is the source. The integral in \autoref{eq: general integral problem green} for $G$ being the Green's function for 3D Laplace is the electric potential at the target point $\x$ due to charges on the source $\Delta$. 
The polynomial factor $p$ is included since most functions can be approximated with a polynomial. 

In \autoref{eq: general integral problem gradient green}, the denominator behaves like $r^3$ as $\x$ and $\y$ becomes close. 
One power is canceled out from the $(\x-\y)$ in the numerator, another is canceled out from $dS(\y)$. When $\y$ approaches $\x$, the vector $(\x - \y)$ becomes orthogonal to $\nx$, so the dot product also removes one power. Thus, the original integrand on the $C^2$ boundary only has a removable singularity. For a more in-depth analysis see chapter 9 from \cite{atkinson1997numerical}. 
Once the smooth boundary is descretized into triangular elements, $\x$ being very close to $\y$ can mean that they are on the same triangle $\Delta$. This results in $(\x-\y)$ being on the same plane as the triangle for all $\y\in\Delta$. If $\Delta$ is not orthogonal to $\nx$, a new singularity is introduced. If $\Delta$ is orthogonal to $\nx$, then the entire integral is zero as the numerator is zero for all $\y\in\Delta$. Both cases result in very bad approximations of the true integral. For a simple visualization of this phenomenon, see \autoref{fig:normaldiagram}.
Even if $\x$ does not lie on $\Delta$, just by being close enough, the integral can numerically appear singular, making the error large.

\begin{figure}
\begin{subfigure}[t]{0.6\textwidth}
    \begin{center}

\begin{tikzpicture}
\begin{axis}[width=8cm,
    height=10cm,
    view={-20}{20},
    colormap/blackwhite,
    axis line style={draw=none},
    tick style={draw=none},
    xticklabel=\empty,yticklabel=\empty,zticklabel=\empty
]
\addplot3[surf,shader=interp,domain=-0.5:1.5,opacity=0.5]{0.5*x*x+0.5*y*y};

\def\xOne{0.0} \def\yOne{0.0}
\def\xTwo{1.0} \def\yTwo{-1.0}
\def\xThree{1.0} \def\yThree{1.0}
\pgfmathsetmacro{\zOne}{0.5*\xOne*\xOne+0.5*\yOne*\yOne}
\pgfmathsetmacro{\zTwo}{0.5*\xTwo*\xTwo+0.5*\yTwo*\yTwo}
\pgfmathsetmacro{\zThree}{0.5*\xThree*\xThree+0.5*\yThree*\yThree}
\pgfmathsetmacro{\dx}{\xOne} 
\pgfmathsetmacro{\dy}{\yOne} 
\def\scalen{-0.5}

\addplot3[only marks, mark=*, mark size=2pt, mark options={fill=red},point meta=explicit symbolic,nodes near coords] coordinates {(\xOne,\yOne,\zOne)[$\x_1$]};
\addplot3[only marks, mark=*, mark size=2pt, mark options={fill=red},point meta=explicit symbolic,nodes near coords]  coordinates {(\xTwo,\yTwo,\zTwo)[$\x_2$]};
\addplot3[mark=*, mark size=2pt, mark options={fill=red},point meta=explicit symbolic,nodes near coords] coordinates {(\xThree,\yThree,\zThree)[$\x_3$]};

\def\yOneX{0.9} \def\yOneY{0.4}
\pgfmathsetmacro{\yOneZ}{0.5*\yOneX*\yOneX+0.5*\yOneY*\yOneY}
\def\yTwoX{0.35} \def\yTwoY{-0.2}
\pgfmathsetmacro{\yTwoZ}{0.5*\yTwoX*\yTwoX+0.5*\yTwoY*\yTwoY}
\def\yHatOneX{0.9} \def\yHatOneY{0.4} \def\yHatOneZ{0.9}
\def\yHatTwoX{0.35} \def\yHatTwoY{-0.2} \def\yHatTwoZ{0.35}

\addplot3[only marks, mark=*, mark size=2pt, mark options={fill=red},point meta=explicit symbolic,nodes near coords] coordinates {(\yOneX,\yOneY,\yOneZ)[$\quad\y_1$]};
\addplot3[only marks, mark=*, mark size=2pt, mark options={fill=red},point meta=explicit symbolic,nodes near coords] coordinates {(\yTwoX,\yTwoY,\yTwoZ)[$\quad\y_2$]};
\addplot3[only marks, mark=*, mark size=2pt, mark options={fill=red},point meta=explicit symbolic,nodes near coords] coordinates {(\yHatOneX,\yHatOneY,\yHatOneZ)[$\quad\Hat{\y}_1$]};
\addplot3[only marks, mark=*, mark size=2pt, mark options={fill=red},point meta=explicit symbolic,nodes near coords] coordinates {(\yHatTwoX,\yHatTwoY,\yHatTwoZ)[$\quad\Hat{\y}_2$]};

\addplot3[patch,shader=interp,mesh/color input=explicit,data cs=cart,opacity=0.5] coordinates {
        (\xOne,\yOne,\zOne) [color=red]
        (\xTwo,\yTwo,\zTwo) [color=red]
        (\xThree,\yThree,\zThree) [color=red]
    };

\pgfmathsetmacro{\norm}{sqrt(\dx*\dx + \dy*\dy + 1)}
\pgfmathsetmacro{\nx}{-\dx/\norm}
\pgfmathsetmacro{\ny}{-\dy/\norm}
\pgfmathsetmacro{\nz}{1/\norm}

\pgfmathsetmacro{\vOneX}{\xTwo-\xOne}
\pgfmathsetmacro{\vOneY}{\yTwo-\yOne}
\pgfmathsetmacro{\vOneZ}{\zTwo-\zOne}
\pgfmathsetmacro{\vTwoX}{\xThree-\xOne}
\pgfmathsetmacro{\vTwoY}{\yThree-\yOne}
\pgfmathsetmacro{\vTwoZ}{\zThree-\zOne}
\pgfmathsetmacro{\dTx}{\vOneY*\vTwoZ - \vOneZ*\vTwoY}
\pgfmathsetmacro{\dTy}{\vOneZ*\vTwoX - \vOneX*\vTwoZ}
\pgfmathsetmacro{\dTz}{\vOneX*\vTwoY - \vOneY*\vTwoX}
\pgfmathsetmacro{\normT}{sqrt(\dTx*\dTx + \dTy*\dTy + \dTz*\dTz)}
\pgfmathsetmacro{\nTx}{\dTx/\normT}
\pgfmathsetmacro{\nTy}{\dTy/\normT}
\pgfmathsetmacro{\nTz}{\dTz/\normT}

\addplot3[-stealth,black,very thick,
    ] coordinates {
        (\xOne,\yOne,\zOne)
        ({\xOne + \nx*\scalen},{\yOne + \ny*\scalen},{\zOne + \nz*\scalen}) 
    };
    
\addplot3[only marks, mark=none, mark size=2pt, mark options={fill=none},point meta=explicit symbolic,nodes near coords] coordinates {({\xOne + \nx*\scalen},{\yOne + \ny*\scalen},{\zOne + \nz*\scalen})[$~~\qquad \boldsymbol{n}(\boldsymbol{x}_1)$]};

\addplot3[-stealth,blue,
    ] coordinates {(\xOne,\yOne,\zOne) (\yOneX,\yOneY,\yOneZ)};
\addplot3[-stealth,blue,
    ] coordinates {(\xOne,\yOne,\zOne) (\yTwoX,\yTwoY,\yTwoZ)};
\addplot3[-stealth,green,
    ] coordinates {(\xOne,\yOne,\zOne) (\yHatOneX,\yHatOneY,\yHatOneZ)};
\addplot3[-stealth,green,
    ] coordinates {(\xOne,\yOne,\zOne) (\yHatTwoX,\yHatTwoY,\yHatTwoZ)};
    
\end{axis} 
\end{tikzpicture}
\end{center}
\vspace{-0.5em}
\subcaption{3D view}
\end{subfigure}%
\hfill
\raisebox{2cm}{\begin{subfigure}[t]{0.3\textwidth}
\begin{center}
\begin{tikzpicture}

\coordinate (bottom) at (0,-2);  
\coordinate (right) at (1.6641, -1.1094); 
\coordinate (middle) at (0.93312114,-1.7689785);
\coordinate (bottom2) at (0, -2.5);
\coordinate (middle2) at (0.9, -1.52);

\draw[thick] (bottom) arc[start angle=-90, end angle=0, radius=2cm];

\draw[fill=red, draw=black, thick] (bottom) circle (1.5pt); 
\node[left] at (bottom) {$\x_1$};

\draw[fill=red, draw=black, thick] (right) circle (1.5pt); 
\node[right] at (right) {$\x_3$};

\draw[fill=red, draw=black, thick] (middle) circle (1.5pt); 
\node[below right] at (middle) {$\y_1$};

\draw[fill=red, draw=black, thick] (middle2) circle (1.5pt); 
\node[above left] at (middle2) {$\hat{\y}_1$};

\draw[-, thick, red] (right) -- (bottom);
\draw[->, thick, black] (bottom) -- (bottom2);
\draw[->, thick, blue] (middle) -- (bottom);
\end{tikzpicture}
\end{center}
\vspace{1.825cm}
\subcaption{2D view}
\end{subfigure}}
\caption{An example showing the issue with approximating a smooth surface with a flat element. (a) The true smooth boundary between the points $\x_1,\x_2,\x_3$ is being approximated with a flat triangular element. On the true smooth boundary, as $\y$ approaches $\x_1$, the vectors $\y - \x_1$ becomes orthogonal to the normal vector $\boldsymbol{n}(\boldsymbol{x}_1)$. However, this is not true when we use the flat element. For any $\Hat{\y}$ on the flat element, the dot product between the vector $\Hat{\y} - \x_1$ and vector $\boldsymbol{n}(\boldsymbol{x}_1)$ is non-zero. 
(b) A 2D side view.
\label{fig:normaldiagram}}
\end{figure}

When evaluating \autoref{eq: general integral problem green} when $\x$ is very close to $\Delta$, previous papers have presented methods involving the Duffy transform \cite{duffy1982quadrature}, polar coordinate systems \cite{cai2002singularity}, and singularity extraction methods \cite{jarvenpaa2003singularity}. Though they work well for integrals of the form \autoref{eq: general integral problem green}, they do not generalize to the case of \autoref{eq: general integral problem gradient green}. This is especially true when $\x$ is in $\Delta$ or on the boundary of $\Delta$. 
If one were to use the recursive formulas of \cite{jarvenpaa2003singularity} for \autoref{eq: general integral problem gradient green} with $\x$ being a vertex of $\Delta$, two terms in their equation (45) will equal to zero. Specifically, they use the divergence theorem to change the integral over $\Delta$ into an integral over the three edges of the triangle. Since the normal vector of these edges are orthogonal to the edge, contributions from two of the edges will be zero. This is undesirable as they are setting something that is non-zero to zero. 

The fundamental problem of the integral \autoref{eq: general integral problem gradient green} when $\x$ is in $\Delta$ is that we are trying to approximate the true weakly singular integral $I_{\text{true}}$ with a strongly singular integral $I = \infty$. Previous analytic methods \cite{jarvenpaa2003singularity} then approximate this divergent integral by computing only its finite part. This results in the loss of information and decreases accuracy. This can also be thought of as an attempt to approximate ``infinity" with a finite value that hopefully is close to $I_{\text{true}}$. 
Instead of approximating this ``infinity" that was intended to approximate $I_{\text{true}}$, we introduce a new evaluation method, which we call Quadratic Surface Approximation. This novel method directly approximates $I_{\text{true}}$ by using geometric information of the true domain $\partial\Omega$. 

The paper is divided into 4 parts. First, \autoref{sec: background} provides some inspiration of why these integrals matter and how they can be applied in boundary element methods (BEM). In \autoref{sec: near singular integrals}, we introduce analytical formulas that can integrate \autoref{eq: general integral problem gradient green} when $\x$ is near $\Delta$ but not on it and the algorithm to determine integration bounds of triangles. In \autoref{sec: singular integrals}, we present Quadratic Surface Approximation. In \autoref{sec: results} we give numerical simulations of our new methods and see how they compare to the standard adaptive quadrature method of MATLAB and setting singular integrals to zero. Finally, we provide an appendix that gives explicit formulas of some integrals that arise in our methods for low degrees of $p$.

An implementation of our new methods in Rust is available on \href{https://github.com/AfZheng126/AnalyticIntegrals}{Github}, and instructions on how to call it using C or Matlab are provided in the repository. 

\section{Background - applications in BEM}\label{sec: background}

It is common practice that partial differential equations (PDEs) are reformulated into integral equations and solved using boundary element methods (BEM). Given a linear PDE boundary value problem, 
\begin{align}
    \mathcal{L}[u] &= 0, \quad \x\in\Omega,\\
    \mathcal{B}[u] &= b(\x), \quad \x\in\partial\Omega,
\end{align}
we can represent the solution $u$ as a single layer potential 
\begin{equation}
    u(\x) = \int_{\partial\Omega} G(\x, \y) \gamma(\y) \diff S(\y),
\end{equation}
in which $G$ is the Green's function corresponding to the differential operator $\mathcal{L}$ and $\gamma$ is some unknown function defined on $\partial\Omega$.
If the boundary condition is Dirichlet $\mathcal{B}[u] = u$, we have integral equations of the form 
\begin{equation}
    \int_{\partial\Omega} G(\x, \y) \gamma(\y)\diff S(\y) = b(\x).
\end{equation}
This first kind integral equation is smoothing, and therefore it is not a good formulation for a BEM. 

If instead the boundary condition is Neumann $\mathcal{B}[u] = \frac{\partial u}{\partial \n}$, then we have integral equations of the form
\begin{equation}
    -\frac{1}{2}\gamma(\x) + \int_{\partial\Omega} \frac{\partial G}{\partial \n(\x)}(\x, \y)\gamma(\y)\diff S(\y) = b(\x).
\end{equation}
This second kind integral equation is much more well-suited for a BEM. 

The boundary $\partial\Omega$ is then descretized into triangles $\Delta$ and nodes $\x_1, \ldots, \x_N$ and it is imposed that the integral equations are satisfied at each $\x_i$.
As $\gamma$ is unknown, it is typical to approximate it with the first $k$ terms of a polynomial basis.
Hence, the integral equations become
\begin{equation}
    -\frac{1}{2} \gamma(\x_j) + \sum_{\Delta} \int_{\Delta} \frac{\partial G}{\partial \n(\x_j)}(\x_j, \y)\left(\sum_{i=1}^k \gamma(v_{\Delta, i} ) l_{\Delta, i}(\y)\right)\diff S(\y) = b(\x_j),
\end{equation}
in which $v_{\Delta, i}$ are some quadrature points in $\Delta$ and $l_{\Delta, i}$ are some polynomial basis functions on $\Delta$. 
This is a linear system of the form $A\gamma = b$, where $A$ is a matrix with entries equal to integrals of the form \autoref{eq: general integral problem green} and \autoref{eq: general integral problem gradient green}. 
To obtain accurate values for $\gamma$, it is vital that all integrals in the matrix $A$ are computed accurately.  

Most common elliptical PDEs have Green's function that contain a $\frac{1}{|\x - \y|}$ singularity. Hence, when we use single layer potentials to solve the PDE with Neumann boundary conditions, we have to consider terms that contain $\frac{(\x - \y)\cdot \n(\x)}{|\x - \y|^3}$. We now present some examples of such Green's functions in $\mathbb{R}^3$.

The Green's function for Helmoltz in $\mathbb{R}^3$ is 
\begin{equation}
    G_H(\x, \y) = \frac{e^{ik|\x - \y|}}{4\pi |\x - \y|},
\end{equation}
and its normal derivative in $\x$ is 
\begin{equation}
    K_H(\x, \y) = \left(-\frac{e^{ik|\x - \y|}}{4\pi |\x - \y|^3} + \frac{ike^{ik|\x - \y|}}{4\pi |\x - \y|^2}\right)(\x - \y) \cdot \n(\x).
\end{equation}
Similar to the Green's function for Laplace, the first term in $K_H$ is singular on flat triangles. 
The reason is also due to the term $(\x - \y)\cdot n(\x)$. On flat triangles, this dot product may not necessarily become zero even when $\y\in\Delta$ approaches $\x$.
We note that our method is effective for low frequencies which is when the polynomial approximation of the exponential term is accurate. The difficulty of high frequency is well known in the literature \cite{babuska1997pollution}. 
The screened Poisson equation's Green's function is 
\begin{equation}
    G(\x, \y) = \frac{e^{-\mu |\x - \y|}}{4\pi |\x -\y|}
\end{equation}
so the integrals are very similar to the Helmoltz case.

Another example is the Cauchy-Navier equations from classical elastostatics \cite{mclean2000strongly, rizzo1983some}. 
The differential operator in this case is 
\begin{equation}
    \mathcal{L}[\boldsymbol{u}] = -\mu\nabla\cdot (\nabla \boldsymbol{u}) - (\mu + \lambda) \nabla(\nabla \cdot \boldsymbol{u})
\end{equation}
in which $\mu$ and $\lambda$ are the Lam\'e coefficients. 
The Green's function for this in $\mathbb{R}^3$ is 
\begin{equation}
    G(\x, \y) = \frac{1}{8\pi\mu (2\mu + \lambda)}\left((3\mu + \lambda) \frac{1}{|\x - \y|} I_3 + (\mu + \lambda) \frac{(\x - \y) (\x - \y)^T}{|\x - \y|^3}\right).
\end{equation}
This Green's function is a $3\times 3$ matrix. Nonetheless evaluating $\int_{\partial\Omega} G(\x, \y) \gamma(\y)\diff S(\y)$ in which $\gamma(\y)$ is now a vector in $\mathbb{R}^3$ involves components of the form \autoref{eq: general integral problem green}. If we consider its normal derivative in $\x$, then we will also see components of the form \autoref{eq: general integral problem gradient green}. 
The PDE governing Stoke's flow is a linear elliptic PDE, and its Green's function is similar to the Green's function for the Cauchy-Navier equation. Hence, our work can also be used to solve Stoke's flow.

\section{Closed form expression for near singular integrals}\label{sec: near singular integrals}

We present a novel method to analytically evaluate our integrals on triangles that works for all cases except for when the target point $\x$ lies inside the triangle $\Delta$ when integrating \autoref{eq: general integral problem gradient green}. 
As we primarily deal with the Green's function of 3D Laplace, we are trying to evaluate 
\begin{equation}\label{eq: Grad Green kernal integral}
    I = \int_{\Delta} \frac{(\x - \y)\cdot \n(\x)}{|\x - \y|^3} p(\y) \diff S(\y), \quad \x \notin \Delta,
\end{equation}
in which $p$ is any polynomial. 
The case when $\x$ lies inside the triangle is left to \autoref{sec: singular integrals} because one of our analytic formulas is divergent. 

Evaluation is aided using a polar coordinate transform that greatly simplifies the integrand. Once in the new polar coordinate system, we give an procedure that finds the bounds of integration on a triangle to obtain explicit closed form formulas for the integral. 

As $\x$ becomes closer to $\Delta$, the numerical integral starts to diverge, making it ``near-singular''. Adaptive refinement of the triangles can give very accurate results in the near-singular case \cite{atkinson1997numerical}, but also greatly increases algorithm runtime. 

Other papers have also given analytic formulas for these integrals over triangles, but either their formulas require special assumptions or need to split the triangle into special sub-triangles \cite{carley2013analytical, okon1982potentialLinear, bohm2024efficient}. Graglia \cite{graglia1993numerical} gives very similar polar transforms as we do for near singular integrals involving Green's functions of Helmholtz and Laplace in 3D, but our new method for determining bounds on integration can be used for the integration of an arbitrary function on triangles. 


\subsection{Choosing the correct polar coordinates}
To obtain an analytic solution to our integrals over a triangle, we first normalize the region of integration via rigid motions in $\mathbb{R}^3$ and then introduce suitable coordinates. 
We map $\Delta$ to a triangle $\Delta_{XY}$ that lies on the $XY$-plane, such that the target point $\x$ lies on the $Z$-axis. 
This map is the rotational map between the triangle normal and the vector $[0, 0, 1]^T$, followed by a translation. 
Let $\n(\Delta)$ be the normal vector of the arbitrary triangle. We want to find a map that rotates $\n(\Delta)$ into $[0, 0, 1]^T$, the normal vector of a triangle that lies on the $XY$-plane. 
This rotational map $R$ is the $3\times 3$ matrix defined using Rodrigues' rotation formula.
After applying the rotation, a translation is done to map the triangle onto the $XY$-plane and the point $R(\x)$ onto the $Z$-axis. 

As these transformations preserve distance, we can now without loss of generality assume that $\x = [0, 0, c]^T$ and the arbitrary triangle $\Delta$ lies on the $XY$-plane. 
The new coordinate system then gives $|\x - \y|^2 = y_x^2 + y_y^2 + c^2$, so \autoref{eq: Grad Green kernal integral} is 
\begin{equation*}
    I = \int_{\Delta_{XY}} \frac{-y_xn_x-y_yn_y+cn_z}{(y_x^2+y_y^2+c^2)^{3/2}} p(\y)\, \diff S_{\y}.
\end{equation*}
If the polynomial $p$ has degree at most 1, we simply need to compute the integrals  
\begin{align}\label{eq: list of near-singular integrals}
I_0:=\int_\Delta \frac{1}{(y_x^2+y_y^2+c^2)^{3/2}} ~dS_{\y},~ I_x:=\int_\Delta \frac{y_x}{(y_x^2+y_y^2+c^2)^{3/2}}~\diff S_{\y}, \\ 
I_y:=\int_\Delta \frac{y_y}{(y_x^2+y_y^2+c^2)^{3/2}}~\diff S_{\y}, I_{x^2}:=\int_\Delta \frac{y_x^2}{(y_x^2+y_y^2+c^2)^{3/2}}~\diff S_{\y}, \\
I_{y^2}:=\int_\Delta \frac{y_y^2}{(y_x^2+y_y^2+c^2)^{3/2}}~\diff S_{\y},~ I_{xy}:=\int_\Delta \frac{y_xy_y}{(y_x^2+y_y^2+c^2)^{3/2}}~\diff S_{\y}.
\end{align}
For example, if we consider $p(\y) = 1$, then \autoref{eq: Grad Green kernal integral} is 
\begin{align*}
I =  cn_z I_0 - n_x I_x - n_y I_y.
\end{align*}
For higher degree polynomials, we would have terms $I_{x^3}, I_{x^2y}, \ldots$, whose equations are more complicated and are not presented here but are listed in the appendix. 

Introducing polar coordinates for $\y$: $y_x=r\cos\theta,~ y_y=r\sin\theta$, we can rewrite the integrals into the form
\begin{equation*}
    \int_{r_\mathrm{start}}^{r_\mathrm{end}} \int_{\theta_\mathrm{start}(r)}^{\theta_\mathrm{end}(r)} \frac{r^{a+b+1}(\cos(\theta))^a(\sin(\theta))^b}{(r^2+c^2)^{3/2}} ~\diff \theta \diff r, 
\end{equation*}
for powers $(a,b)$, the powers of $y_x$ and $y_y$ in the numerator, and $r_\mathrm{start}$, $r_\mathrm{end}$, $\theta_\mathrm{start}(r)$, $\theta_\mathrm{end}(r)$ the bounds of integration that represent $\Delta$. 
The method to determine the values of $\theta_\mathrm{start}(r)$ and $\theta_\mathrm{end}(r)$ is explained in the following subsection.

\subsection{Determining start and end angles}\label{sec: finding bounds of integration}
When we integrate in polar coordinates starting at the point $\y$ where the circle of radius $r$ intersects an edge of the triangle, 
the angle $\theta(r)$ of the point $\y$ can be written as 
\begin{equation}
    \theta(r) = \arccos\left(\frac{y_x}{r}\right).
\end{equation}
However, $y_x$ is a non-linear function of $r$, so this does not give rise to a simple analytic formula. 
Instead, we look at the projection of the origin to the lines which define the sides of the triangle. 
Let us label the positively oriented vertices as $\V_1, \V_2, \V_3$. 
Let $\boldsymbol{d}_1$ be the projection of the origin to the line on which $\V_1$ and $\V_2$ lie on; $\boldsymbol{d}_2$ be the projection of the origin to the line on which $\V_2$ and $\V_3$ lie on; 
and $\boldsymbol{d}_3$ be the projection of the origin to the line on which $\V_3$ and $\V_1$ lie on.

Let $\phi_k$ be the angle between $\boldsymbol{d}_k$ and the $X$-axis, then the angle at a point of intersection between a circle of radius $r$ and a side of the triangle is 

\begin{equation}\label{equation for boundary theta}
    \theta(r) = \pm\arccos\left(\frac{d_k}{r}\right) + \phi_k, 
\end{equation}
in which $d_k$ is the norm of the orthogonal projection of the origin onto the line on which that specific edge of the triangle lies on. An example is shown in \autoref{fig:regiondiagram}. 
The importance of this equation is that now 
\begin{equation*}
    \theta_\mathrm{start}(r) = \mathrm{sign}_\mathrm{start}\arccos\left(\frac{d_\mathrm{start}}{r}\right) + \phi_\mathrm{start}, 
\end{equation*}
in which $\mathrm{sign}_\mathrm{start}, d_\mathrm{start}, \phi_\mathrm{start}$ are all piece-wise constant in $r$. 
If we split our integral into the correct regions in $r$, it is much simpler to write out an analytic solution to our integral. 

\begin{figure}
\begin{center}
\begin{tikzpicture}

\begin{scope}[xshift = -4cm, yshift=-2cm]
\coordinate (xaxis) at (4.5,0);
\coordinate (yaxis) at (0,3.5);
\draw[->] (0,0) -- (xaxis) node[right] {$x$};
\draw[->] (0,0) -- (yaxis) node[above] {$y$};

\coordinate (V1) at (1.0,1.5);
\coordinate (V2) at (4.0,2.5);
\coordinate (V3) at (1.5,3.0);
\coordinate (dP12) at ( 3.0, 1.0);
\coordinate (dP23) at (-2.5, 0.5);
\coordinate (dP31) at (-0.5,-1.5);
\coordinate (P12r) at (-2,0.5);
\draw[fill=blue!20,opacity=0.5] (V1) -- (V2) -- (V3) -- cycle;

\node at (V1) [above left] {$\boldsymbol{V}_1$};
\node at (V2) [above] {$\boldsymbol{V}_2$};
\node at (V3) [above] {$\boldsymbol{V}_3$};
\fill (V1) circle (2pt);
\fill (V2) circle (2pt);
\fill (V3) circle (2pt);

\def\lval{3.0}; 
\coordinate (L12) at (2.2953,1.9318);
\coordinate (L13) at (1.3867,2.6602);
\fill (L12) circle (1pt);
\fill (L13) circle (1pt);
\node at (L12) [left] {};
\node at (L13) [left] {};

\coordinate (O) at (0,0);
\fill (O) circle (2pt);
\draw[thick] (L12) arc (40.0847:62.4676:\lval);    
\draw[] (O) -- ++(0.6256*\lval,0.7802*\lval) node[below] {$r$};
\draw[] (O) -- (L12);
\draw[] (O) -- (L13);

\coordinate (Q12) at (-0.3500, 1.0500);
\coordinate (Q23) at ( 0.6346, 3.1731);
\coordinate (Q31) at ( 0.4500,-0.1500);

\fill (Q12) circle (2pt);
\draw[dashed] (O) -- ++(Q12) node[midway, left] {$d_1$};
\draw[dashed] (P12r) -- (Q12) -- (V1);
\pic [draw,below right, ->, "$~\phi_{1}$", angle eccentricity=1.0, angle radius=15] {angle = xaxis--O--Q12};
\pic [draw,right, ->, "$\theta_{12}(r)$", angle eccentricity=1.0, angle radius=50] {angle = xaxis--O--L12};

\pic [draw,above left, <-, "$\arccos(d_{1}/r)~~$", angle eccentricity=1.2, angle radius=28] {angle = L12--O--Q12};
\end{scope}

\begin{scope}[xshift = 4cm]

\coordinate (xaxis) at (3.2,0);
\coordinate (yaxis) at (0,3.2);
\draw[->] (-3.2,0) -- (xaxis) node[right] {$x$};
\draw[->] (0,-3.2) -- (yaxis) node[above] {$y$};

\coordinate (V1) at (0, 2);
\coordinate (V2) at (-2, -1);
\coordinate (V3) at (3, 0);
\coordinate (dP12) at (-0.9231, 0.6154);
\coordinate (dP23) at (0.1154, -0.577);
\coordinate (dP31) at (0.923, 1.3846);

\def\rone{0.5883}; 
\def\rtwo{1.1094}; 
\def\rthree{1.6641}; 
\def\rfour{2}; 
\def\rfive{2.2361}; 
\def\rsix{3}; 

\draw[fill=blue!20, opacity=0.5] (0, 0) circle (\rsix);
\draw[fill=orange!20, opacity=0.5] (0, 0) circle (\rfive);
\draw[fill=blue!20, opacity=0.5] (0, 0) circle (\rfour);
\draw[fill=orange!20, opacity=0.5] (0, 0) circle (\rthree);
\draw[fill=blue!20, opacity=0.5] (0, 0) circle (\rtwo);
\draw[fill=orange!20, opacity=0.5] (0, 0) circle (\rone);

\node at (V1) [above left] {$\boldsymbol{V}_1$};
\node at (V2) [above] {$\boldsymbol{V}_2$};
\node at (V3) [above] {$\boldsymbol{V}_3$};
\node at (dP12) [above] {$\boldsymbol{d}_1$};
\node at (dP23) [above] {$\boldsymbol{d}_2$};
\node at (dP31) [above] {$\boldsymbol{d}_3$};
\fill (V1) circle (2pt);
\fill (V2) circle (2pt);
\fill (V3) circle (2pt);
\fill (dP12) circle (2pt);
\fill (dP23) circle (2pt);
\fill (dP31) circle (2pt);
\draw[] (V1) -- (V2) -- (V3) -- cycle;
\end{scope}

\end{tikzpicture}
\end{center}
\caption{On the left, we show how to define the variables for the integration. The triangle $\Delta$ is defined by the points $\boldsymbol{V}_1, \boldsymbol{V}_2, \boldsymbol{V}_3$. 
The distance from the origin to the projection of the origin onto the side $\overline{\boldsymbol{V}_1\boldsymbol{V}_2}$ is $d_{1}$; this projection makes an angle of $\phi_{1}$ with the $X$-axis. 
A circle of radius $r$ centered at the origin intersects the side $\overline{\boldsymbol{V}_1\boldsymbol{V}_2}$ with an angle of $\theta_{12}(r)$ with the $X$-axis. The three labeled angles give that $\theta_{12}(r) = \phi_1 - \arccos(d_1/r)$. 
On the right, we show the different critical radii for integration. 
The projection of the origin to the edges of the triangles are labeled with $\boldsymbol{d}_1, \boldsymbol{d}_2, \boldsymbol{d}_3$, which in this case all lie on their respective edges. 
Letting $-1$ represent inactive, $1$ represent active, and $0$ represent split, the activities of the edges ordered in $[\overline{\boldsymbol{V}_1\boldsymbol{V}_2}, \overline{\boldsymbol{V}_2\boldsymbol{V}_3}, \overline{\boldsymbol{V}_3\boldsymbol{V}_1}]$ change as follows. 
Starting from $r=0$, $[-1, -1, -1]\rightarrow [-1, 0, -1] \rightarrow [0, 0, -1] \rightarrow [0, 0, 0] \rightarrow [1, 0, 1] \rightarrow [-1, 1, 1] \rightarrow [-1, -1, -1]$.}\label{fig:regiondiagram}
\end{figure}

The critical radii $R_i$ which separates the regions on which $\mathrm{sign}_\mathrm{start}, d_\mathrm{start}, \phi_\mathrm{start}$ (and the corresponding end values) are constant are the norm of the three vertices of the triangle 
and the norm of the orthogonal projections $\boldsymbol{d}_k$ if they lie on the triangle, ordered from smallest to largest. 
Starting at $r = 0$, the $\theta$ integral is either on $[0, 2\pi]$ or on $\emptyset$, depending on if the origin is in the triangle. 
This is easily checked numerically. Now as we progress through the different critical radii, we need to keep track of which angles we start and stop the integration. 
This is done by labeling the triangle edges and keeping track of their ``activity", which is defined as follows. 
\begin{definition}
    For a given radius $r$, the edge $\overline{\boldsymbol{V}_1\boldsymbol{V}_2}$ of a triangle with vertices $\boldsymbol{V}_1, \boldsymbol{V}_2, \boldsymbol{V}_3$ lying on the $XY$-plane is \textbf{inactive} if the circle of radius $r$ on the $XY$ plane does not intersect $\overline{\boldsymbol{V}_1\boldsymbol{V}_2}$. 
    It is \textbf{active} if the circle of radius $r$ intersects it only once, and \textbf{split} if it intersects it twice. 
\end{definition}
The reason for this definition is to know which angles we use in our integrals. 
For example, at $r=0$, all edges are inactive unless the origin is one of the vertices, in which case two edges are active and one is inactive. 
If two edges are active while the third is inactive, then we only need to consider two angles and which one of them is the start angle and which is the end angle. 

To make sure that we are always integrating in the correct direction, we need to first order the vertices of the triangle such that its vertices $\boldsymbol{V}_1, \boldsymbol{V}_2, \boldsymbol{V}_3$ are positively oriented and $\boldsymbol{V}_1$ is the vertex closest to the origin. 

\begin{definition}\label{def: left and right}
    For a triangle in the previously specified configuration, if the projection point $\boldsymbol{d}_1$ lies on the interior of the edge $\overline{\boldsymbol{V}_1\boldsymbol{V}_2}$, it splits the edge into $\overline{\boldsymbol{V}_1\boldsymbol{d}_1}$ and $\overline{\boldsymbol{d}_1\boldsymbol{V}_2}$, 
    which we define as its \textbf{left} and \textbf{right} respectively. For other projection points, the left and right are defined so that everything is still positively oriented. 
    If the projection point does not lie on the edge, that entire edge is considered as both left and right for the sake of the method. 
\end{definition}

Note that when a projection point lies on an edge of the triangle, its left and right corresponds to different signs for arc-cosine. 
It is not trivial to determine which has the positive sign and which has the negative sign, so it must be computed by testing which sign of arc-cosine the vertices of that edge use. 
If $\boldsymbol{d}_1$ lies on $\overline{\V_1\V_2}$, the simplest way is to check if which sign makes $\V_1 = \pm\arccos\left(\frac{d_1}{\|\V_1\|}\right) + \phi_1$ true.
This configuration makes it so that we always integrate from edges $\overline{\V_1\V_2}$ to $\overline{\V_2\V_3}$ to $\overline{\V_3\V_1}$. 
In the case where an edge is split, we use \autoref{def: left and right} to know that we always integrate starting from the right side of the first active or split edge and then alternate between left and right of the following edges. 

As the radius $r$ increases, the activity of the edges change when $r$ equals to a critical radii. 
At $r=R_i$, the change in the activities is dependent on whether $R_i$ is the norm of a vertex or a projection point. 
If it is the norm of a vertex, the activities of the edges the vertex lie on changes. Inactive edges become active while active edges become inactive. 
A split edge on the other hand becomes active. There is a special case where two vertices have the exact same norm, which can cause a split edge to become inactive. 
However, due to floating point precision, we do not consider this case. If $R_i$ is the norm of a projection point, the only case that happens is the inactive edge that it lies on becomes split. 
An example of how the activities change is shown in \autoref{fig:regiondiagram}. 

This method can then be used to integrate any function on any triangle that lies on $\mathbb{R}^2$ without the need to apply some non-distance preserving transformation. Rust code that implements this is provided in \href{https://github.com/AfZheng126/AnalyticIntegrals}{Github} with documentation.

\subsection{Example computation of an integral}

Now that we can accurately determine the activity of each edge and which sign each arc-cosine has, we can integrate on the triangle. We denote
\begin{equation*}
\begin{split}
    \theta_{\mathrm{end}}(r) &= \mathrm{sign}_\mathrm{end}\arccos\left(\frac{d_\mathrm{end}}{r}\right) + \phi_\mathrm{end},\\
    \theta_{\mathrm{start}}(r) &= \mathrm{sign}_\mathrm{start}\arccos\left(\frac{d_\mathrm{start}}{r}\right) + \phi_\mathrm{start}.
\end{split}
\end{equation*}

The analytic formula for $I_0$ is 
\begin{equation}\label{eq: I0}
\begin{split}
    I_0  & = \int_{r_\mathrm{start}}^{r_\mathrm{end}} \int_{\phi_\mathrm{start} + \mathrm{sign}_\mathrm{start}\arccos\left(d_\mathrm{start}/r\right)}^{\phi_\mathrm{end} + \mathrm{sign}_\mathrm{end}\arccos\left(d_\mathrm{end}/r\right)} \frac{r}{(r^2 + c^2)^{3/2}}\, \diff \theta \, \diff r \\
    & = (\phi_\mathrm{start} - \phi_\mathrm{end}) \int_{r_\mathrm{start}}^{r_\mathrm{end}} \frac{r}{(r^2+c^2)^{3/2}}\, \diff r + \mathrm{sign}_\mathrm{end}\int_{r_\mathrm{start}}^{r_\mathrm{end}} \frac{r}{(r^2 + c^2)^{3/2}}\arccos\left(\frac{d_\mathrm{end}}{r}\right)\, \diff r \\
    & - \mathrm{sign}_\mathrm{start}\int_{r_\mathrm{start}}^{r_\mathrm{end}} \frac{r}{(r^2 + c^2)^{3/2}}\arccos\left(\frac{d_\mathrm{start}}{r}\right)\, \diff r.\\
\end{split}
\end{equation}
The three integrals in $r$ have analytic formulas which are found in \autoref{appendix: formula for integral of normal derivative of single layer potential}. The other integrals of $I_x, I_y , \ldots$ once written in a similar form can also be evaluated by analytically evaluating the integrals in $r$ by hand or by using symbolic software like Mathematica. 

Similar formulas can also be obtained for \autoref{eq: general integral problem green} integrals and can work for any arbitrary $\x$.
As $c$, $d_\mathrm{start}$ and $d_\mathrm{end}$ become extremely small, numerical errors can occur in terms that involve $\arctan$ or $\log$. In this case, we use the first four terms of the Taylor expansion. 

\subsection{Why this fails for targets on the triangle}
Let us consider $I_0$ when $\x\in\Delta$.
Since $\x$ is translated to the origin, $\Delta$ contains the origin. 
Hence, $r_\mathrm{start} = c = 0$.
If we look at the first integral in \autoref{eq: I0}, we get  
\begin{equation*}
    \int_{0}^{r_\mathrm{end}} \frac{1}{r^2} ~\diff \theta \diff r. 
\end{equation*}
This is undefined and shows why this method only works for $\x\notin\Delta$.


\section{Quadratic Surface Approximation}\label{sec: singular integrals}
In the previous section, we were able to directly integrate \autoref{eq: general integral problem gradient green} when $\x\notin\Delta$ because the integrand was only singular at $\y = \x\notin\Delta$. Once $\x\in\Delta$, the integral becomes strongly singular and all the previous formulas fail. Specifically, one of the integrals in $r$ become divergent. As the pre-discretization integral is clearly not infinite \cite{atkinson1997numerical}, we stress that this is due to the process of approximating a $C^2$ boundary with flat triangles, and not a feature of the integral we are trying to approximate. 

Let $\V_1, \V_2, \V_3$ be three points on the surface $\partial\Omega$ and let $\partial\Omega_V\subset\partial\Omega$ be the curved triangular patch on which we want to integrate our function. The specifics of how to define this triangular patch are discussed later. The true integral we want to evaluate is 
\begin{equation}\label{eq: True grad of Green integral on C2 boundary}
    I_{\text{true}} = \int_{\partial\Omega_V} \frac{\partial G}{\partial n(\x)}(\x, \y) p(\y)\, \diff S(\y).
\end{equation}
In this section, we provide a novel method to approximate $I_{\text{true}}$ when $\x$ is on $\partial\Omega_V$ or near it. 
To calculate this integral, we need to either approximate the domain of integration, or approximate the integrand. 

If we approximate $\partial\Omega_V$ with a flat triangle $\Delta$, this integral becomes divergent. Thus, we need to integrate a new function on $\Delta$ whose integral value is finite and hopefully approximates $I_{\text{true}}$. Some simply integrate the zero function on $\Delta$, or only evaluate the finite part of the divergent integral \cite{jarvenpaa2003singularity}. This approach can be viewed as trying to approximate infinity with a finite value that hopefully is close to $I_{\text{true}}$.
Previous methods that keep the true domain $\partial\Omega_V$ do not have explicit generalizable formulas \cite{mantivc1994computing}, or are applicable only in the $2$D setting with line elements \cite{sladek2000optimal}. 

We first present Quadratic Surface Approximation, a method which keeps the domain to be $\partial\Omega_V$ and approximates the integral with analytic closed form expressions. Then we give analytic proofs for why it works. Finally we conclude with how to numerically compute certain geometric data of $\partial\Omega$ needed for Quadratic Surface Approximation. 

\subsection{Motivation}
Instead of approximating an ``infinity" that was intended to approximate $I_{\text{true}}$, Quadratic Surface Approximation directly approximates $I_{\text{true}}$ by using geometric properties of the true domain $\partial\Omega$. 

Let $\V_1, \V_2, \V_3$ be three points on the surface $\partial\Omega$ and $\partial\Omega_V\subset\partial\Omega$ be the curved triangular patch bounded by the geodesics between the three points. 
Suppose the target point $\x$ is sufficiently close by to $\partial\Omega_V$. We wish to compute 

\begin{equation}
    I_{\text{true}} = \frac{1}{4\pi}\int_{\partial\Omega_V} \frac{(\x - \y)\cdot \n(\x)}{\|\x - \y\|^3} \gamma(\y) \diff S(\y).
\end{equation}

The idea is to locally represent the boundary $\partial\Omega$ locally as a two-dimensional manifold defined as the graph of a quadratic function. 
This allows us to approximate many terms in the integral using Taylor series while also incorporating geometric properties of the manifold. 
However, the way to do this can be tricky, so we split it into a two cases.

We first consider the case where $\x$ is in $\partial\Omega$ and close to $\partial\Omega_V$. 
This means that $\n(\x)$ is defined with respect to $\Omega$. It is fine if $\x$ does not lie on the triangular patch $\partial\Omega_V$, as long as it is close to it. 
Then we consider the case when $\x\notin \partial\Omega$, so that $\n(\x)$ may just be some vector that is not related to $\Omega$. 

\subsection{Case 1: When the target point is on the boundary}

We first rotate and translate the coordinate system so that $\x$ becomes the origin and its normal vector $\n(\x)$ becomes $[0, 0, 1]^T$.
Now the $XY$-plane is tangent to $\partial\Omega$, $\partial\Omega$ can be locally represented as the graph of a function $f(s_1, s_2)$. 
This is where we require that $\x$ is close to $\partial\Omega_V$, otherwise this $f$ may not be an accurate representation of the manifold at $\partial\Omega_V$. 
These rotations do not guarantee that the flat triangle $\Delta$ lies on the $XY$-plane,  
but this is not a concern as we are trying to integrate on $\partial\Omega_V$ and not $\Delta$. 
Let $p(\V_1\V_2\V_3)$ be the projection of $\partial\Omega_V$ onto the $XY$-plane. 
Note that there is no reason why $p(\V_1\V_2\V_3)$ is a triangle, as it is just a projection of geodesics. 
Let us denote $p\Delta$ as the triangle with vertices $p\V_i$, the projections of $\V_i$ onto the $XY$-plane.
The notation here shows that $p\Delta$ is just the projection of $\Delta$ to the $XY$-plane. 

The true integral over $\partial\Omega_V$ can now be written as an integral on $p(\V_1\V_2\V_3)$ using a pushforward map. 
\begin{equation}
    I_{\text{true}} =\frac{1}{4\pi}\int_{p(\V_1\V_2\V_3)} \pi_*\left[\frac{(\x - \y)\cdot \n(\x)}{\|\x - \y\|^3} \gamma(\y)\right] \pi_*[\diff S(\y)].
\end{equation}

When evaluating, we first note that $p(\V_1\V_2\V_3) \ne p\Delta$. 
Let 
\begin{equation}
    \epsilon = \max_{(s_1, s_2)\in p\Delta} \sqrt{s_1^2 + s_2^2},
\end{equation}
so that $p\Delta$ is contained in some $B_\epsilon$ ball centered around the origin.
We have that
\begin{align*}
    &\left|\int_{p(\V_1\V_2\V_3)} \pi_*[K(\x, \y)\gamma(\y)] \pi_*[\diff S(\y)] - \int_{p\Delta} \pi_*[K(\x, \y)\gamma(\y)] \pi_*[\diff S(\y)]\right|\\
    &\quad \leq \int |\chi_{p(\V_1\V_2\V_3)}(\y) -\chi_{p\Delta}(\y)||\pi_*[K(\x, \y)\gamma(\y)]| \pi_*[\diff S(\y)].
\end{align*}

Assuming that the manifold is smooth, the closer the triangles are to the origin, the flatter they are and the more they agree. 
As the triangles become finer or as $\|\y\|$ decreases, this error will also decrease with a very rough error bound proportional to the area of the triangular patch which is $\mathcal{O}(\epsilon^2)$.
However, this error is accounted for when integrating over the entire boundary. 
When we create a triangular mesh, we are plotting many points on the boundary of the surface, and then connecting them with curves.
On individual triangles there will be a difference in domain of integration, but that difference will be accounted for in the neighbouring triangles. 
As the rotations we perform are with respect to $\x$, rotations will not affect this mis-alignment of triangular domains which are in $\y$.
Thus it is fine to integrate on $p\Delta$ and not $p(\V_1\V_2\V_3)$ when considering the integration over the entire boundary, otherwise the error is roughly bounded by $\mathcal{O}(\epsilon^2)$.

Next, we claim is that we can approximate $\pi_*[\diff S(\y)]$ with $\diff S(\y)$ with $\mathcal{O}(\epsilon^2)$ error.
This will be very useful as $\diff S(\y)$ is much more easier to compute than its pushforward. 
To see why this is true, we introduce some manifold theory. 
Given an $n$-dimensional manifold $\mathcal{M}$ embedded in $\mathbb{R}^m$ parametrized locally by $h(s_1, \ldots, s_n)\in\mathbb{R}^m$,
the partial derivatives of $h$ are vectors of length $m$. 
The first fundamental form can be represented by the symmetric matrix defined by 
\begin{equation}
    g_{ij} = \frac{\partial h}{\partial s_i} \cdot \frac{\partial h}{\partial s_j}
\end{equation}
This then gives the volume form on $\mathcal{M}$ as $\diff S_{\mathcal{M}}(\y) = \sqrt{\det g} \, \diff s_1\wedge\cdots \wedge \diff s_n$.

In the case where the manifold is the graph of a function in $\mathbb{R}^{n+1}$, we have that $h(s_1, \ldots, s_n) = (s_1, \ldots, s_n, f(s_1, \ldots, s_n))$.
Calculations then give 
\begin{align*}
    \diff S_{\mathcal{M}}(\y) &= \sqrt{\det g}\, \diff s_1\wedge\cdots \wedge \diff s_n \\
    &= \sqrt{\det (I + \nabla f \otimes \nabla f)} \, \diff s_1\wedge\cdots \wedge \diff s_n\\
    &= \sqrt{1 + \|\nabla f\|^2}\, \diff s_1\wedge\cdots \wedge \diff s_n.
\end{align*}
In our case, $\partial\Omega$ is parameterized as the graph of a function, so 
\begin{equation}
    \pi_*[\diff S(\y)] = \sqrt{1 + \|\nabla f\|^2}\, \diff s_1\wedge\cdots \wedge \diff s_n.
\end{equation}
Let $m(\y) = \sqrt{1 + \|\nabla f(\y)\|^2}$. 
By the positioning of the coordinate system, we have that $f(\underline{0}) = 0, \nabla f(\underline{0}) = \underline{0}$, so the Taylor expansion of $m$ is 
\begin{align*}
    m(\y) &= m(\underline{0}) + \sum_{i=1}^{n} \frac{\partial m}{\partial s^i}(\underline{0}) y_i + \cdots\\
    &= 1 + \cdots
\end{align*}
Hence, $1 - \mathcal{O}(\epsilon^2) \leq m(\y) \leq 1 + \mathcal{O}(\epsilon^2)$ and $\pi_*[\diff S(\y)] = \diff S(\y) + \mathcal{O}(\epsilon^2)$.
Finally, we claim that 
\begin{equation}
    \begin{split}
        &\frac{1}{4\pi}\int_{p\Delta} \pi_*\left[\frac{(\x - \y)\cdot \n(\x)}{\|\x - \y\|^3} \gamma(\y)\right] \diff S(\y)\\
        &\quad =\frac{1}{4\pi}\int_{p\Delta} \frac{\left(- \frac{1}{2}\sum_{i, j = 1}^{2}K_{ij}s_is_j\right)\pi_*\left[\gamma(\y)\right]}{\|\x - \y\|^3} \diff S(\y) + \mathcal{O}(\epsilon^2),
    \end{split}
\end{equation}
in which $K_{ij}$ are the second order derivatives of $f$, but we write them as $K_{ij}$ as they are also the second fundamental forms of the manifold $\partial\Omega$ at $\x$. 

Intuitively, the second fundamental form is an extrinsic curvature that describes how curved a surface is. 
For a more rigorous and geometric understanding, see \cite{toponogov2006differential, do2016differential} or any textbook on Differential Geometry. 
Given a point $\x\in\partial\Omega$ such that locally, $\partial\Omega$ can be represented as the graph of the function $f(s_1, s_2)$ the second fundamental form at a point $(s_1, s_2, f(s_1, s_2))$ is the bilinear form defined as
\begin{equation}\label{eq: second fundamental form for graph}
    II(u, v) = u\cdot \frac{1}{\sqrt{1 + f_{s_1}^2 + f_{s_2}^2}}\begin{pmatrix}
        f_{s_1 s_1} & f_{s_1 s_2} \\
        f_{s_2 s_1} & f_{s_2 s_2}
    \end{pmatrix} v.
\end{equation}
With our rotations, $\x$ is now the origin, and the first order partial derivatives of $f$ are zero at the origin, so the middle term on the right hand side becomes the hessian of $f$. 
Thus, the second order derivatives of $f(s_1, s_2)$ are equivalent with the second fundamental form evaluated at $\x$ locally.
Intuitively, when the manifold lies tangent to the $(s_1, s_2)$ plane, its hessian describes the local curvature of the manifold. 
The reason for the introduction of the second fundamental form is for ease of computation.
BEM makes the domain of integration a two dimensional manifold, but computationally all the coordinate systems are in three dimensions.
Hence, finding the function $f$ can be hard, and taking its derivatives even more so.
The graph representation of the manifold allows us to skip the step of finding the graph function $f$ and just compute the second fundamental form of the manifold. 
Furthermore, the second fundamental form can be defined using the shape operator for easier calculations. 

Given a point $\x\in\partial\Omega$, the shape operator is simply $\mathcal{S}(\x) = -\nabla \n(\x)$, which is just the gradient of the normal vector. 
When our boundary is given by a pseudo-signed distance function $F(x, y, z)$, the normal vector at any boundary point is just the normalized gradient of $F$, so 
\begin{equation}\label{eq: shape operator}
    \mathcal{S}(\x) = - \frac{\nabla^2 F(\x)}{\|\nabla F\|}.
\end{equation}
The second fundamental form is defined as 
\begin{equation}
    II(u, v)(x) = u\cdot \mathcal{S}(\x) v
\end{equation}
This new definition can be shown to be equivalent to \autoref{eq: second fundamental form for graph} by explicitly calculating $\n(\x)$ and taking its gradient.

When computing $K_{ij}$, we simply choose $u, v$ to be $e_i, e_j$, the standard basis in which we integrate. 
However, it is important to know that the basis in which we integrate is not the same as the original standard basis vectors in $\mathbb{R}^3$.
This is because we had to rotate our coordinate system. 
Therefore, the actual second fundamental form is computed via
\begin{equation}\label{eq: second fundamental form}
    \Pi(u, v)(x) = (R_1^{-1}\cdots R_m^{-1} u)\cdot \mathcal{S}(\x) (R_1^{-1}\cdots R_m^{-1} v)
\end{equation}
where $R_i$ were the rotations we performed to make $\n(\x)$ be $[0, 0, 1]^T$. 
Translations do not influence the basis functions in which we integrate in, so there is no need to consider those. 
Instead of having to find a function $f$ whose graph locally represents our manifold and taking its second order derivatives, 
it may be easier to calculate the shape operator as it only uses information of the normal vectors $\n(\x)$. 

As our triangulation method approximates $\gamma$ with linear functions, we shall now regard $\gamma$ as a polynomial $p$. 
Using our previous interpolation method, we have that 
\begin{equation}
    \pi_*[p](x^1, x^2) \approx \sum_{k=1}^M \pi_*[p](pQ_k)l_k(x^1, x^2),
\end{equation}
for some quadrature points $pQ_k\in p\Delta$. Each $pQ_k$ is the projection of a point $Q_k\in \partial\Omega_V$ onto the $XY$-plane.
By definition of the push-forward, we get that this is equivalent to
\begin{equation}
    \pi_*[p](s_1, s_2) \approx \sum_{k=1}^M p(Q_k)l_k(s_1, s_2).
\end{equation}
If we choose a quadrature with enough nodes, this is exact for our polynomial $p$. 
For each $k$, we need to show that 
\begin{equation*}
    \left|\pi_*\left[\frac{(\x - \y)\cdot \n(\x)}{\|\x - \y\|^3}\right] l_k(s_1, s_2) \diff S(\y) - \frac{\left(- \frac{1}{2}\sum_{i, j = 1}^{2}K_{ij}s_is_j\right)}{\|\x - \y\|^3}l_k(s_1, s_2) \diff S(\y)\right| \leq C\epsilon^2
\end{equation*}
for some finite $C>0$. 
To do this, we use the fact that $\x$ is the origin, $\n(\x)$ is $[0, 0, 1]$, and $\y$ can be represented locally as the graph of a function $f(s_1, s_2)$. 
First, we bring the absolute value inside the integrals and bound the integral above by a circle $B$ that contains $p\Delta$. 
This way, we naturally change to the polar coordinates of $(s_1, s_2) = (r\cos(\theta), r\sin(\theta))$. 
Now let $\epsilon$ be the radius of the bounding circle, which is also equal to $\max_{\y\in p\Delta} \|\y - \x\|$. 
The polynomials $l_k$ are smooth functions, so they obtain a maximum inside $B$ and we simply need to bound 
\begin{equation}\label{eq: epsilon error integral}
    \int_B \left| \pi_*\left[\frac{(\x - \y)\cdot \n(\x)}{\|\x - \y\|^3}\right] - \frac{\left(- \frac{1}{2}\sum_{i, j = 1}^{2}K_{ij}s_is_j\right)}{\|\x - \y\|^3}\right| \diff S(\y).
\end{equation}
After explicitly writing out the pushforward map, the term inside the integral is 
\begin{equation}
    \frac{-f(s_1, s_2)}{\left(s_1^2 + s_2^2 + f(s_1, s_2)^2 \right)^{3/2}} - \frac{-\frac{1}{2} \sum_{i,j=1}^{2} \frac{\partial^2 f}{\partial s_i \partial s_j} s_is_j}{(s_1^2 + s_2^2)^{3/2}}.
\end{equation}
Since $f(0, 0) = 0$ and its first order partial derivatives are zero, we can Taylor expand $f$ as 
\begin{equation}
    f(s_1, s_2) = \frac{1}{2} \sum_{i,j=1}^{2} \frac{\partial^2 f}{\partial s_i \partial s_j} s_i s_j + \mathcal{O}(r^3) 
\end{equation}
in which $r = \sqrt{s_1^2 + s_2^2}$. 
For notation, let $Q = \frac{1}{2} \sum_{i,j=1}^{2} \frac{\partial^2 f}{\partial s_i \partial s_j} s_i s_j$ so that $f = Q + \mathcal{O}(r^3)$. 
We also have
\begin{equation}
    \left(s_1^2 + s_2^2 + f^2\right)^{3/2} = \left(r^2 + f^2\right)^{3/2} = r^3 \left(1 + \frac{f^2}{r^2}\right)^{3/2}.
\end{equation}
Since $f = \mathcal{O}(r^2)$, we conclude that 
\begin{equation}
    \left(s_1^2 + s_2^2 + f^2\right)^{3/2} = r^3(1 + \mathcal{O}(r^2)).
\end{equation}
Hence, 
\begin{equation}
    \left(s_1^2 + s_2^2 + f^2\right)^{3/2} - \left(s_1^2 + s_2^2\right)^{3/2} = r^3(1 + \mathcal{O}(r^2) - 1) = \mathcal{O}(r^5). 
\end{equation}
We now simplify the fractions 
\begin{equation}\label{eq: epsilon error 1}
    \begin{split}
        &\frac{-f(s_1, s_2)}{\left(s_1^2 + s_2^2 + f(s_1, s_2)^2 \right)^{3/2}} - \frac{-\frac{1}{2} \sum_{i,j=1}^{2} \frac{\partial^2 f}{\partial s_i \partial s_j} s_is_j}{(s_1^2 + s_2^2)^{3/2}}\\
        &\quad = \frac{-Q + \mathcal{O}(r^3)}{r^3 + \mathcal{O}(r^5)} - \frac{-Q}{r^3}.
    \end{split}
\end{equation}
For the first fraction, we can further simplify by
\begin{equation}\label{eq: epsilon error 2}
    \begin{split}
        \frac{Q + \mathcal{O}(r^3)}{r^3 + \mathcal{O}(r^5)} &= \frac{Q}{r^3}\left(1 + \mathcal{O}(r^2)\right)^{-1} + \frac{\mathcal{O}(r^3)}{r^3}\left(1 + \mathcal{O}(r^2)\right)^{-1}\\
        &= \frac{Q}{r^3}(1 - \mathcal{O}(r^2)) + \mathcal{O}(1) (1 - \mathcal{O}(r^2))\\
        &= \frac{Q}{r^3} (1 - \mathcal{O}(r^2)) + \mathcal{O}(1) + \mathcal{O}(r^2).         
    \end{split}
\end{equation}
The equality on the second line is due to the Taylor expansion of $\frac{1}{1+x}$ for $x$ near zero. 
Plugging \autoref{eq: epsilon error 2} into \autoref{eq: epsilon error 1}, we obtain 
\begin{equation}
    \begin{split}
        \frac{-Q + \mathcal{O}(r^3)}{r^3 + \mathcal{O}(r^5)} - \frac{-Q}{r^3}&= -\frac{Q}{r^3} (1 - \mathcal{O}(r^2)) + \mathcal{O}(1) + \mathcal{O}(r^2) - \frac{-Q}{r^3}\\
        &= -\frac{Q}{r^3}\left(1 - \mathcal{O}(r^2) + 1\right) + \mathcal{O}(1) + \mathcal{O}(r^2)\\
        &= \frac{Q}{r^3}\mathcal{O}(r^2) + \mathcal{O}(1) + \mathcal{O}(r^2)\\
        &= \mathcal{O}(r) + \mathcal{O}(1) + \mathcal{O}(r^2).
    \end{split}
\end{equation}
Using this, \autoref{eq: epsilon error integral} is 
\begin{equation}
    \begin{split}
        &\int_B \mathcal{O}(r) + \mathcal{O}(1) + \mathcal{O}(r^2) \diff S(\y)\\
        &\quad =\int_{0}^{2\pi} \int_{0}^{\epsilon} \left(\mathcal{O}(r) + \mathcal{O}(1) + \mathcal{O}(r^2)\right)r \diff r \diff \theta\\
        &\quad =\mathcal{O}(\epsilon^2).
    \end{split}
\end{equation}

Hence, we conclude that 
\begin{equation}
    I_{\text{true}} = \frac{1}{4\pi}\int_{p\Delta} \frac{- \sum_{i, j = 1}^{2}K_{ij}s_is_j}{(s_1^2 + s_2^2)^{\frac{3}{2}}} \pi_*\left[\gamma(\y)\right] \diff S(s_1)\diff S(s_2) + \mathcal{O}(\epsilon^2)
\end{equation}
where $\pi_*\left[\gamma(\y)\right]$ is just a polynomial. In our case, it is a linear function.
Changing to polar coordinates, this gives us integrals of the form
\begin{equation}
    \int_{p\Delta} r^{a+b+1-3}\cos^a(\theta)\sin^b(\theta) \, \diff \theta\, \diff r,
\end{equation} 
for $a, b$ positive integers such that $a+b\geq 2$. 
As $a+b+1 - 3 \geq 0$, there are no potential divergent integrals in $r$ like in \autoref{sec: near singular integrals}. 

The method to determine the start/stop angles and radii are the same as how we did it for the analytic equations.
If $\x \in \partial\Omega_V$, then $p\Delta$ contains the origin. Otherwise, $p\Delta$ does not contain the origin. 

In the special case when $\x$ is a vertex of $\partial\Omega_V$, we can instead integrate in a simpler way with an extra rotation. 
This commonly occurs in numerical algorithms \cite{hsiao2008boundary} as the points $\x$ are usually collocation points chosen to be vertices of the triangles in $\partial\Omega_{\Delta}$.
Instead of integrating in $\theta$ first like before, we integrate in $r$ first and then integrate in $\theta$. 
The reason for this is because $p\V_1$ is the origin. 
Having a fixed vertex at the origin means that we can perform a rotation so that $p\V_2$ lies on the positive $s_1$-line, which then means that the integral in $\theta$ is a single segment unlike in previous cases. 
Let $p\Delta$ be positively oriented and use a rotation such that $p\V_2$ lies on the $s_1$-axis. 
This way, $\theta\in[0, \theta_{p\V_3}]$ where $\theta_{p\V_3}$ is the angle of $p\V_3$ and the $s_1$-axis. 
As a function of $\theta$, we have that 
\begin{equation}
    r(\theta) =  \frac{|p\V_2| \sin(\theta_{2})}{\sin(\theta + \theta_{2})},
\end{equation}
in which $\theta_{2}$ is the angle between $p\V_1p\V_2$ and $p\V_2p\V_3$. This equation is obtained using the sine law, which states that 
\begin{equation}
    \frac{\sin(\theta_{2})}{r} = \frac{\sin(\pi - \theta - \theta_{2})}{|p\V_2|}.
\end{equation}
For a visualization, see \autoref{fig: radius as a function of theta}. 
\begin{figure}
\begin{center}
\begin{tikzpicture}

\coordinate (xaxis) at (4.5,0);
\coordinate (yaxis) at (0,3.5);
\draw[->] (0,0) -- (xaxis) node[right] {$s_1$};
\draw[->] (0,0) -- (yaxis) node[above] {$s_2$};

\coordinate (V1) at (0.0,0.0);
\coordinate (V2) at (4.0,0.0);
\coordinate (V3) at (2.0,3.0);
\coordinate (midP) at (2.2, 2.7);
\draw[fill=blue!20,opacity=0.5] (V1) -- (V2) -- (V3) -- cycle;

\node at (V1) [above left] {$\boldsymbol{p\V_1}$};
\node at (V2) [above right] {$\boldsymbol{p\V_2}$};
\node at (V3) [above] {$\boldsymbol{p\V_3}$};
\fill (V1) circle (2pt);
\fill (V2) circle (2pt);
\fill (V3) circle (2pt);

\node at (midP) [above right] {$\boldsymbol{W}$};
\fill (midP) circle (2pt);
\draw[dashed] (V1) -- (midP);

\draw pic["$\theta_2$", draw=black, <-, angle eccentricity=1.2, angle radius=1cm] {angle = V3--V2--V1};
\draw pic["$\psi$", draw=black, ->, angle eccentricity=1.2, angle radius=0.7cm] {angle = V1--midP--V2};
\draw pic["$\theta$", draw=black, ->, angle eccentricity=1.2, angle radius=0.7cm] {angle = V2--V1--midP};

\end{tikzpicture}
\end{center}
\caption{A diagram showing how to calculate $r(\theta)$. 
As $\theta$ grows, the sine law gives that the ratio between $r$ and $\sin(\theta_2)$ is the same as the ratio between $|p\V_2|$ and $\sin(\psi)$.} \label{fig: radius as a function of theta}
\end{figure}
Hence the integrals we need to compute are of the form
\begin{equation}
    \int_0^{\theta_\mathrm{end}} (\cos(\theta))^a (\sin(\theta))^b \int_0^{r(\theta)} r^{a+b+1-3}\, \diff r \, \diff \theta.
\end{equation}
The integrals in $r$ are trivial as $a+b+1-3\geq 0$, so we get 
\begin{equation}
    \int_0^{\theta_\mathrm{end}} (\cos(\theta))^a (\sin(\theta))^b \frac{1}{a+b-1}\left(\frac{|p\V_2|\sin(\theta_2)}{\sin(\theta+\theta_2)}\right)^{a+b-1}\, \diff \theta.
\end{equation}
The desired integrals can then be written as some linear combination of this integral with $(a, b)\in \{(2, 0), (1, 1), (0, 2), \ldots, (0, 3)\}$. The formulas for these integrals can be found in \autoref{appendix: formula for quadratic approximation}.

\subsection{Case 2: When the target is not on the boundary}
In this case, $\x$ is not in $\partial\Omega_V$ or $\partial\Omega$. 
It is a freely chosen point in $\mathbb{R}^3$, so the normal vector $\n(\x)$ now has nothing to do with $\partial\Omega_V$. 
However, we still assume that $\x$ is relatively close to $\partial\Omega_V$ so that we can locally represented $\partial\Omega$ as the graph of a function $f(s_1, s_2)$. 

As $\x$ is no longer on $\partial\Omega_V$, it does not make sense to change the coordinate system so that $\x$ is the origin.
Instead, we rotate and translate the coordinate system so that $\x$ becomes $[0, 0, c]^T$ and the manifold which is the graph of $f(s_1, s_2)$ is tangent to the $XY$-plane at the origin. 
In this setting, the origin is the projection of $\x$ onto $\partial\Omega$ and $c$ is the distance between $\x$ and $\partial\Omega$.

Previously, the second fundamental form was calculated at $\x$. 
As $\x$ is no longer on $\partial\Omega$, we instead need to compute the second fundamental form at the projection of $\x$ onto $\partial\Omega$. 
As $\n(\x)$ is just some arbitrary constant now, we can disregard it and only look at the vector 
\begin{equation}
    \x - \y = \begin{bmatrix}
        -s_1\\
        -s_2\\
        c - f(s_1, s_2)
    \end{bmatrix}.
\end{equation}
Using the same steps, we still have very similar claims, but we conclude that 
\begin{equation}
    I_{\text{true}} = \frac{1}{4\pi}\int_{p\Delta} \frac{\begin{bmatrix}
    -s_1\\
    -s_2\\
    c - \frac{1}{2} \sum_{i, j = 1}^{2}K_{ij}s_is_j
    \end{bmatrix}\cdot \n(\x)}{(s_1^2 + s_2^2 + c^2)^{\frac{3}{2}}} \pi_*\left[\gamma(\y)\right] \diff S(\y) + \mathcal{O}(\epsilon^2).
\end{equation}
Since $c \ne 0$, the denominator actually never becomes zero, so the proof is actually a lot easier as we can just approximate each term in the integral with its Taylor expansions.
Specifically, letting $Q = \frac{1}{2} \sum_{i,j=1}^{2} K_{ij} s_is_j$, we have 
\begin{equation}
    \begin{split}
        &\left|\frac{c-f}{\left(s_1^2 + s_2^2 + (c-f)^2\right)^{3/2}} - \frac{c - Q}{\left(s_1^2 + s_2^2 + c^2\right)^{3/2}} \right|\\
        & \quad \leq \left| \frac{c-f}{\left(s_1^2 + s_2^2 + (c-f)^2\right)^{3/2}} - \frac{c-Q}{\left(s_1^2 + s_2^2 + (c-f)^2\right)^{3/2}}\right| \\
        & \quad + \left| \frac{c-Q}{\left(s_1^2 + s_2^2 + (c-f)^2\right)^{3/2}} - \frac{c - Q}{\left(s_1^2 + s_2^2 + c^2\right)^{3/2}} \right|.
    \end{split}
\end{equation}
The first term is easily seen to be $\mathcal{O}(r^2)$.
The second term is also $\mathcal{O}(r^2)$ by comparing their Taylor expansions.
The fact that $c\ne 0$, allows the Taylor expansion to be defined locally.

Converting to polar coordinates, we still have integrals of the form
\begin{equation}
    \int_{p\Delta} \frac{r^{a+b+1} \cos^a(\theta)\sin^b(\theta)}{\left(r^2 + c^2\right)^{\frac{3}{2}}}\, \diff \theta\, \diff r,
\end{equation} 
for $a, b$ non-negative integers. 
Since $c \ne 0$, none of the integrals diverge. 
We have already dealt with these types of integrals in the analytic methods, but now we need to consider cases when $a, b =2$.
These are listed in the appendix. 

These ideas of the second fundamental form can also be applied to other Green's functions.
The main idea is to just use rotations to simplify $(\x- \y)$ and $\n(\x)$ and then represent $\y$ as the graph of a function $f(s_1, s_2)$. 
Since the rotations make the manifold tangential to the $XY$-plane, the second order derivatives of $f$ become the second fundamental form of the manifold. 
The integration bounds of the projected triangle $p\Delta$ are then found with the algorithm shown in \autoref{sec: finding bounds of integration}. 

\subsection{How to compute the second fundamental form}
If a signed-distance function is available, then one can first compute the shape operator using \autoref{eq: shape operator}, and then compute the second fundamental form using \autoref{eq: second fundamental form}.

In many cases, a signed-distance function or some kind of function parameterization of the interfaces are not available.
Instead, only a file containing locations of each node and the edges between them are given.
Thus, we need a way to obtain geometric information from a list of nodes (vectors in $\mathbb{R}^3$), and a list of triangles (vectors in $\mathbb{Z}_+^3$ specifying which nodes are the vertices).
Only the shape operator is needed for the second fundamental form, so we just need to compute the partial derivatives of the normal vector on the interface.
To calculate these partial derivatives, we first need to approximate the normal vectors at each node. 
The simplest way to do this approximation is to take some weighted average of the normals of each triangle that the node is a vertex on.
These triangle normals are calculated using the cross product, but it is important to check that the normals are pointing outwards.
The algorithm itself is unable to know which direction is inward or outward, so we have to assume that the list of triangles are given so that the node labels in each vector allow the cross product to point outward. 
Hence,
\begin{equation}
    \n(\x) \approx \frac{\sum_{\Delta: \x\in\Delta} \alpha_\Delta \n(\Delta)}{\| \sum_{\Delta: \x\in\Delta} \alpha_\Delta \n(\Delta)\| }
\end{equation}
A common choice for $\alpha_\Delta$ is the incident angle of $\x$ on $\Delta$. The bigger the angle, the greater the weight. 

Once the normal vectors are approximated, the next step is to approximate the partial derivatives. 
Due to the lack of information of the normal vectors at non-node points, standard finite difference schemes do not work. 
Instead, we need to choose 3 different nodes such that the three edges formed at not collinear and solve
\begin{gather}
    \n(\x_1) \approx \n(\x) + \nabla\n(\x) \cdot (\x_1 - \x), \\
    \n(\x_2) \approx \n(\x) + \nabla\n(\x) \cdot (\x_2 - \x), \\
    \n(\x_3) \approx \n(\x) + \nabla\n(\x) \cdot (\x_3 - \x).
\end{gather}
This can be converted into a $9\times 9$ linear system on the components of $\nabla\n(\x)$. 
Furthermore, with some re-ordering of the rows, we actually get a matrix composed of three 3x3 blocks on the diagonal. 
Two of these points can be chosen to be nodes such that $(\x, \x_1, \x_2)$ is one the the triangles, while the third is on a different triangle that $\x$ is a vertex of. 
As long as all three edges are not collinear, this small linear system can be quickly solved using gaussian elimination. 

A naive implementation of this gives large errors. 
The reason for this is because the triangles are representing a 2D surface, so the third direction is very close to being a linear combination of the other two directions.
Hence, the last three eigenvalues of the matrix is a lot smaller than the other 6.

To deal with this ill-conditioned problem, we need to recognize that Quadratic Surface Approximation does not need information of the gradient in the third direction. 
Though the shape operator describes how the normal vector changes in the entire space, the second fundamental form only needs to know how the normal vector changes along the boundary which is a lower dimension. 

\subsection{Polynomial Method to approximate Second Fundamental Form} \label{sec: approximate second fundamental form using polynomials}
If the shaper operator is only used to calculate the second fundamental form, it is better to directly approximate the second fundamental form instead. 
From \autoref{eq: second fundamental form for graph}, we saw that the second fundamental form of a graph is just the hessian of the function $f$. Hence, all we need to do is to compute the second order derivatives of $f$. If $f$ is locally approximated by a polynomial, we just need the coefficients of its second order terms. 

To approximate the second fundamental form at a node $\x$, we need five neighboring nodes $\x_1, \ldots, \x_5$, but not their normals.
First, we translate the coordinate system so that $\x$ is the origin.
Then we rotate the coordinate system so that $\n(\x)$ becomes $[0, 0, 1]^T$. 
This allows us to fit the polynomial $f(s_1, s_2) = as_1^2 + bs_2^2 + cs_1s_2 + ds_1 + es_2$ to our five translated and rotated points. 
Solving this $5\times 5$ linear system for the coefficients, we conclude that the second fundamental form is 
\begin{equation}
    II(u, v)(\x) = u\cdot \begin{bmatrix}
        2a & c\\
        c & 2b
    \end{bmatrix}v.
\end{equation}
Unlike the previous method, we get a $2\times 2$ matrix as we are already restricting ourselves to the two dimensional manifold. The values of $d, e$ are not used because the manifold being tangent to the plane implies that they are equal to zero. 
Here, $u, v$ are the basis vectors in which we integrate in, which are not necessarily the standard basis vectors. 
This is because $K_{ij}$ are the affected by any rotations we make in the $s_1s_2$-plane. 
The first rotation $R_1$ in 3D to make the plane tangent to the manifold does not need to be considered because
\begin{equation}
    R_1^T \mathcal{S}(\x) R_1 = \begin{bmatrix}
        2a & c & \cdot\\
        c & 2b & \cdot\\
        \cdot & \cdot & \cdot
    \end{bmatrix},
\end{equation}
in which $\cdot$ represents values that we do not care about as they are derivatives in the third direction $\n(\x)$.

There are cases where a node lives only on four triangles that combined only have five distinct vertices (like the top of a pyramid).
Instead of just looking at nodes that are connected to $\x$ by a single edge, we may need to consider nodes that are connected to $\x$ by a path of two edges. 
This also gives us more data, so instead of solving a $5\times 5$ linear system, we can find the least squares solution of an over-determined linear system.
Assuming the triangular mesh is fine enough, the least squares solution will be more stable and more accurate, but will take longer to solve for. Numerical tests show that the least squares system always results in better accuracy without a significant increase in runtime.

\section{Results}\label{sec: results}

We test our new methods based on their accuracy which is the relative absolute difference of the values computed via our method and Matlab's \texttt{integral2} function \cite{shampine2008matlab}, which uses an adaptive quadrature. 
It is not the best method to integrate the integrals we consider, but is widely accepted to be a relatively robust method for non-singular integrals. 
We cannot directly compare runtime because our methods are implemented in Rust, which is much faster than Matlab.
However, it should not be too difficult to believe that evaluating a few closed form equations is much faster than adaptive refinement. 

\subsection{Comparisons on near-singular integrals}
Let us consider evaluating the integral 
\begin{equation}\label{eq: K for simulation}
    I = \int_{\Delta} \frac{(\x-\y)\cdot \nx}{4\pi |\x - \y|^3} \, \diff S(\y),
\end{equation}
in which $\Delta$ is a triangle in the triangulation approximation of the boundary $\partial\Omega = S^1$.
We first present our results on the analytic equations given to evaluate the case when $\x$ does not lie in $\Delta$.
As stated in previous sections, the integral is not singular at all, but its behaviour is near-singular. 
The integral is evaluated in 2000 different tests. 
In each test, the first two components of $\Delta$'s vertices are uniformly sampled from $[-0.01, 0.01]$ while the third component is calculated so that the vertices lie on $S^1$. The point $\x$ is chosen to be $[0, 0, c]^T$ where $c$ is uniformly sampled from $[-0.01, 0.01]$, but its normal vector $\n(\x)$ is fixed to be $[\sqrt{3}/3, \sqrt{3}/3, \sqrt{3}/3]^T$. This is because after rotations, the normal vectors $\n(\x)$ are different in each test case. 
\begin{figure}
    \centering
    \includegraphics[width=1.0\linewidth]{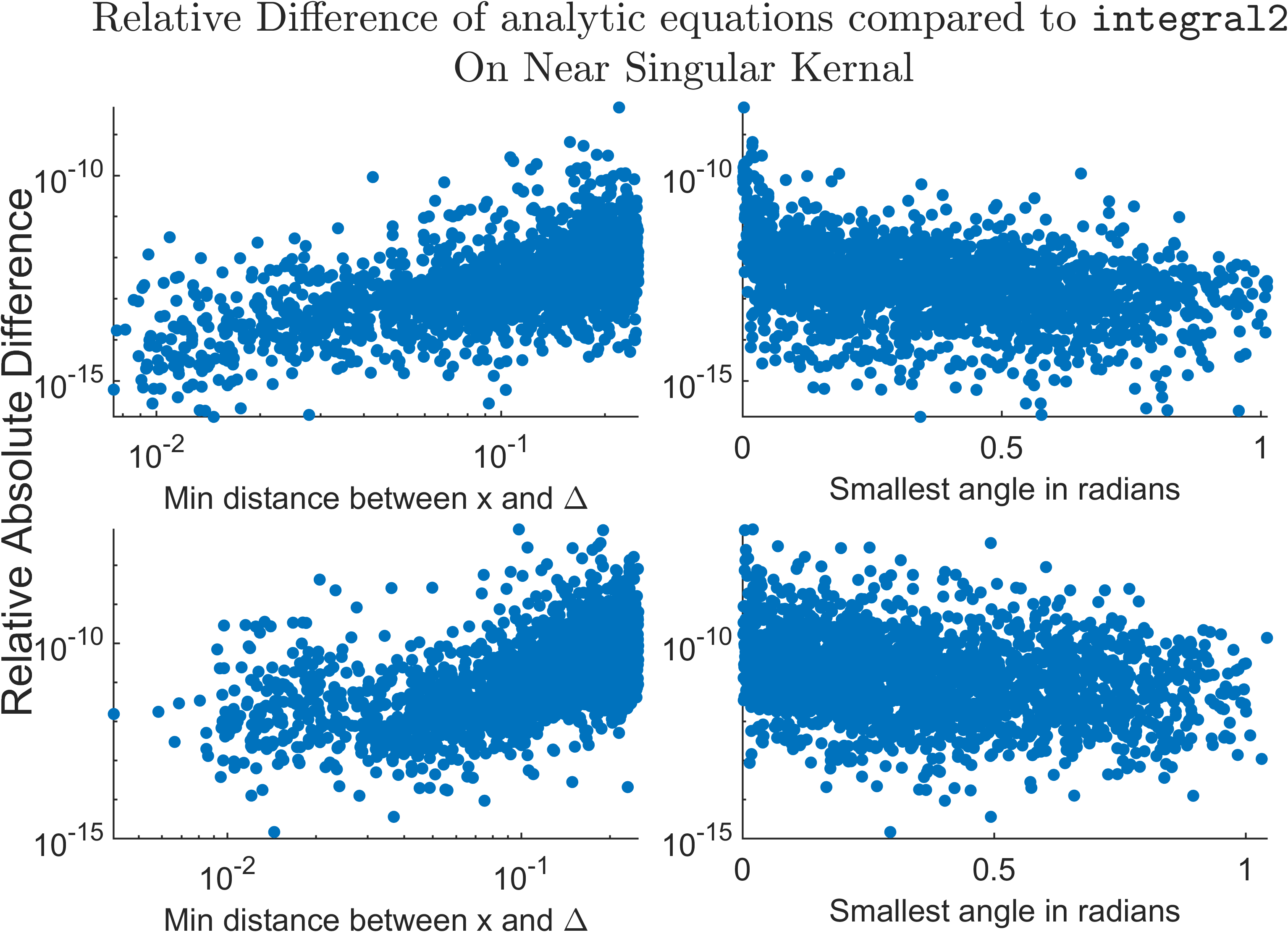}
    \caption{Comparison between the values of integrals using analytic equations and an adaptive method \texttt{integral2} from Matlab on 4000 tests (2000 in each row). The first row plots the results on integrating \autoref{eq: K for simulation}. The second row plots the results on integrating \autoref{eq: K*y for simulation}. $\Delta$ is assumed to have vertices on $S^1$ with their first two components uniformly distributed in $[-0.01, 0.01]^2$. $\x = [0, 0, c]^T$ where $c$ is uniformly distributed in [-0.01, 0.01]. Absolute and relative tolerance for \texttt{integral2} was set to $10^{-14}$.} 
    \label{fig:near singular analytic results}
\end{figure}
The results of these tests are shown in \autoref{fig:near singular analytic results}. 
Even when calculating
\begin{equation}\label{eq: K*y for simulation}
    I = \int_{\Delta} \frac{(\x-\y)\cdot \nx}{4\pi |\x - \y|^3} y_x\, \diff S(\y),
\end{equation}
in which $y_x$ is the first component of $\y$, we obtain very similar results. 

For integrals that do not have strong singularities at all, such as \autoref{eq: general integral problem green}, its similar analytic equations give the same accuracy. Though these equations give values that are very accurate, this analytic method is slower than a simple 3-point quadrature method due to the algorithm to find the bounds of integration.  
It is up to the users to decide when to use this analytic method and when to use a $n$-point quadrature method with small $n$. 

\subsection{Comparisons on singular integrals}
In the next test, we still consider \autoref{eq: K for simulation} but $\x$ is now one of its vertices. 
In each test case, $\Delta$ is simulated the same way as the previous tests, but now $\x$ and $\nx$ are fixed so that $\x$ is the first vertex of $\Delta$ and its normal is the normal vector of $S^1$ at $\x$. 
When running the adaptive method, we actually integrate $K(\x, \y)$ on the curved triangular patch on $S^1$. This is because the integral on $\Delta$ is divergent.

\begin{figure}
    \centering
    \includegraphics[width = 0.9\linewidth]{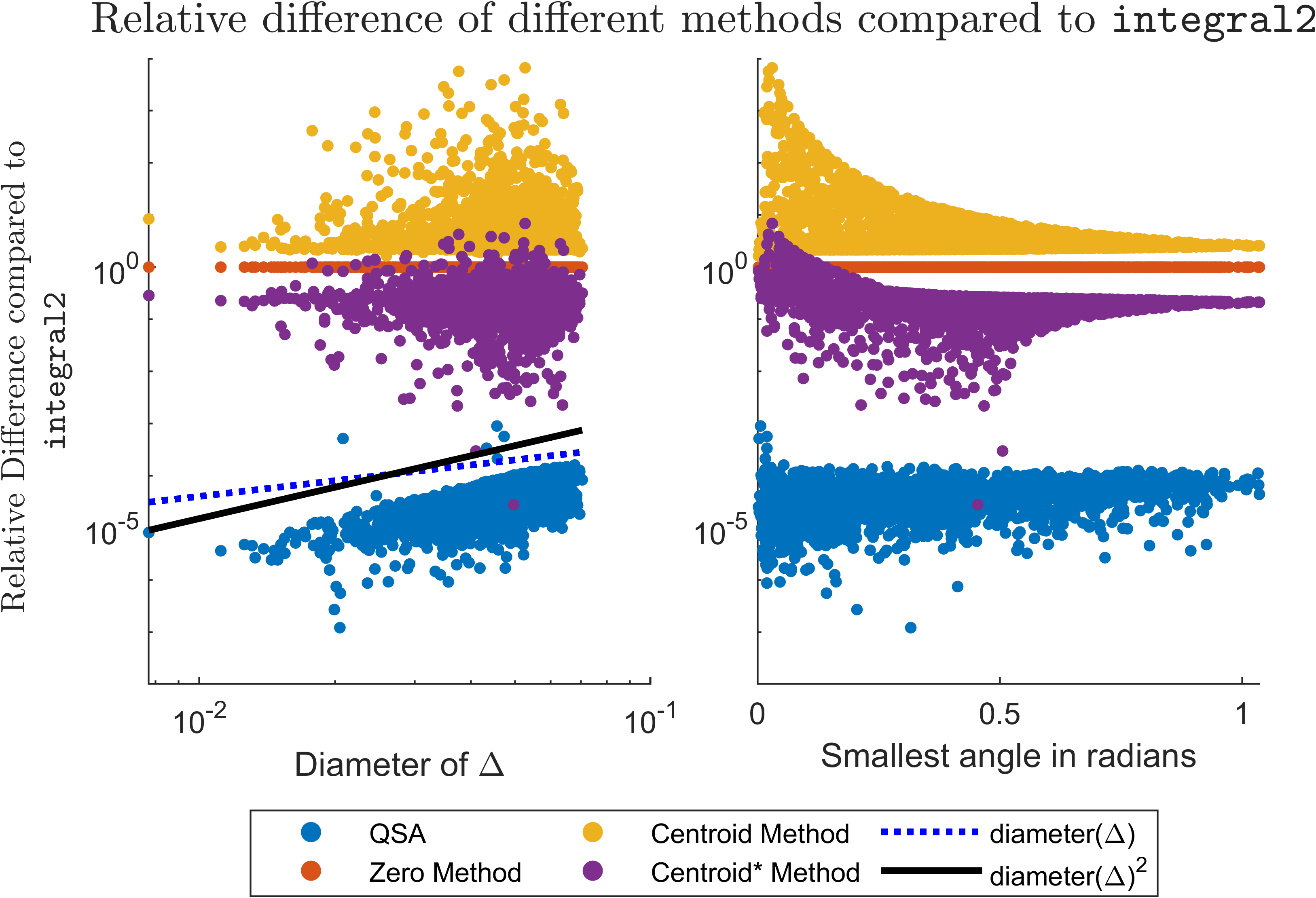}
    
    \caption{Comparison between Quadratic Surface Approximation (QSA), the Zero method (setting the singular integrals to zero), and two Centroid methods that are just 1 point quadrature schemes. The Centroid* method evaluates the function at the projection of the triangle centroid to the true boundary. The Centroid method computes the function at the centroid of the flat triangle. The relative difference of computing the integral \autoref{eq: K for simulation} when $\x$ is a vertex of $\Delta$ on 2000 tests using QSA compared to MATLAB's adaptive method \texttt{integral2}. In the left subplot, the X-axis is the diameter of $\Delta$. The blue dashed line represents a line of slope one, while the black dashed line has slope two. On right subplot, the X-axis is the smallest angle of $\Delta$ in radians. Both subplots show that QSA has a much smaller relative difference compared to the Zero method, though its error depends on the diameter of $\Delta$. The slope of the line that gives an upper bound for the difference for QSA seems to have slope 2. 
    The right graph shows that QSA's accuracy is not dependent on how acute or equilateral the triangles are, but the two centroid methods do depend on it. 
    Absolute and relative tolerance for \texttt{integral2} was set to $10^{-8}$.} 
    \label{fig:singular error result}
\end{figure}

For Quadratic Surface Approximation, the closer the region of integration is to the origin, the better it approximates the adaptive method. This corresponds to the fact we are just using the second order Taylor expansion of $f$, which gives larger errors when $(s_1, s_2)$ becomes further away from the origin.
In the log-log plot, we see that the line that upper bounds the difference has a slope of two, which is what our analysis in \autoref{sec: singular integrals} stated.
The right subplot shows that there is no major trend between difference and smallest angle. 
This suggests that Quadratic Surface Approximation does not depend on how acute or equilateral the triangle is.
The two 1-point quadrature methods do heavily depend on the shape of the triangle. 
The more acute the triangles, the greater the error.
As the Zero method simply sets the integral to zero, its relative difference compared to \texttt{integral2} will always be 1.0. 
We also see that using a 1 point quadrature with the quadrature point being the centroid of the flat triangle actually gives greater differences than the Zero method. 
If the quadrature point was chosen to be the centroid of the curved triangular patch, the difference decreases but can still be worse than the Zero method for extremely acute triangles. 

In general, if any quadrature rule is used, we will encounter the problem of whether the quadrature points are on the true geometry or the flat triangle. Putting the quadrature points on the flat triangle gives approximation error, but putting the quadrature points on the curved triangular patch gives two other errors.
First, if $\x$ coincides with one of the quadrature points, we cannot evaluate the function at that singularity. Secondly, the quadrature weights were designed for flat triangles, so they are wrong for the curved triangular patch. In \autoref{fig:singular error result}, the quadrature weight was always the area of the flat triangle. This is because for such small triangles, the difference in surface area between the flat triangle and the curved triangular patch are on the order of $\mathcal{O}(\epsilon^2)$. Finding the correct quadrature weight for the curved triangular patch is non-trivial. 
Sometimes only a mesh of the true geometry is given, so without a functional parametrization of $\partial\Omega$, the quadrature points can only be on the flat triangle, which gives large approximation error. 
This may explain why the Zero method is commonly used in current literature. 
The results shown in \autoref{fig:singular error result} allow us to conclude that Quadratic Surface Approximation is more accurate than the other methods. 


It is also important to note that if the tolerance is decreased, \texttt{integral2} will start to give warnings. Hence, it is not known if \texttt{integral2} is the best approximation of the true integral values.
Unfortunately, analytic equations for the true solution are not available for triangular patches on $S^1$, so we can only compare to \texttt{integral2}. 

We also consider the case of simply integrating $K_{\mathbb{R}^3}$ on the unit sphere with $\x\in S^2$.
This can be done analytically to get 
\begin{equation}\label{eq: integral of K on sphere}
   \int_{S^2} \frac{(\x - \y) \cdot \n(\x)}{4\pi |\x - \y|} \diff S(\y) = 0.5, \quad \x\in S^2.
\end{equation}
The proof of this equality is found in \autoref{appendix: proof of integrals on the unit sphere}. 
\begin{figure}
    \centering
    \includegraphics[width=0.7\linewidth]{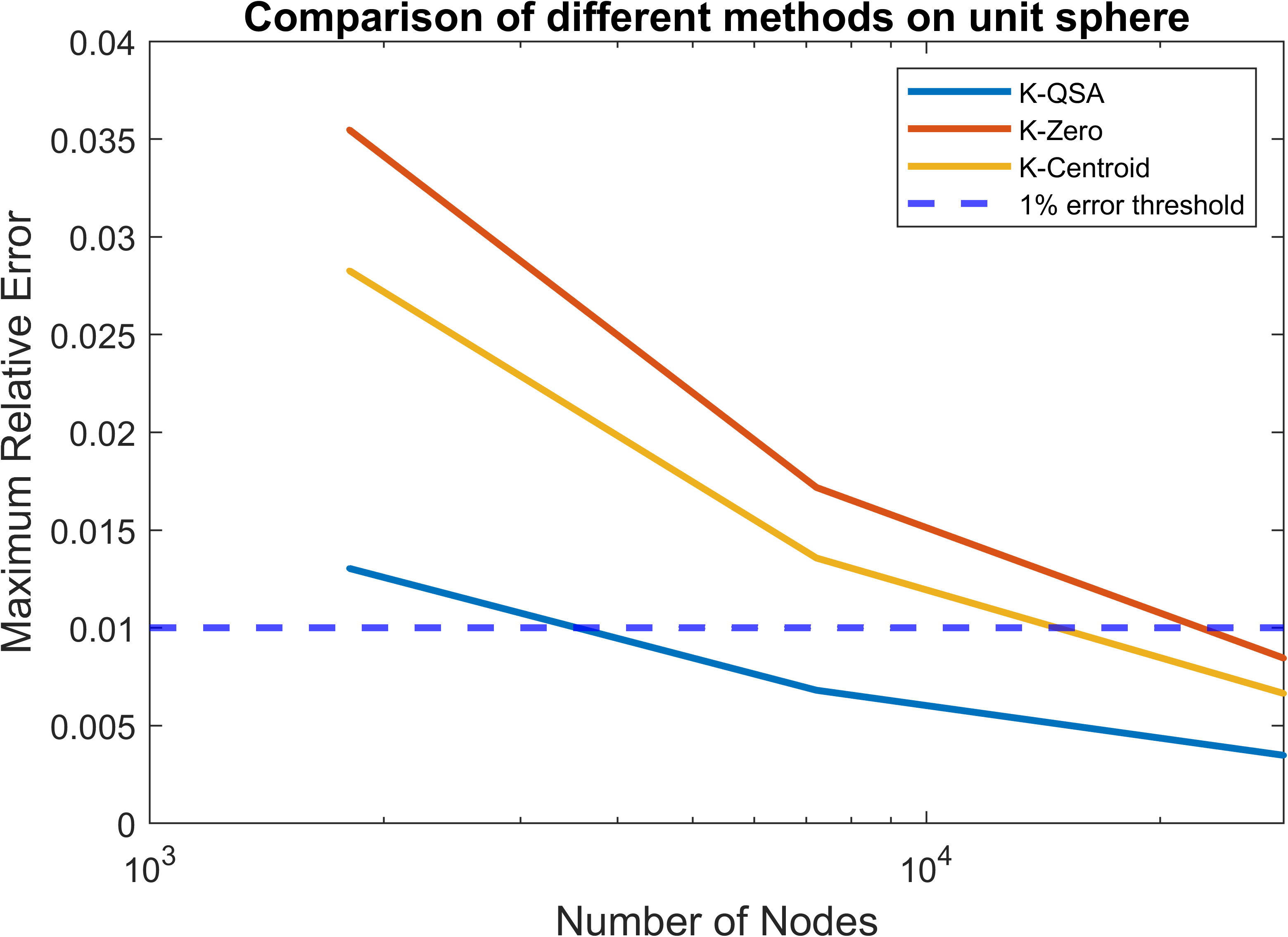}
    \caption{The maximum relative error of integrating $K_{\mathbb{R}^3}$ on the unit sphere at all triangle vertices $\x$ in the triangulation of the unit sphere. All non-singular integrals were evaluated using the analytic equations of \autoref{sec: near singular integrals}. The singular integrals (when $\x$ is a vertex of the triangle $\Delta$) were evaluated using either Quadratic Surface Approximation (QSA), setting it to zero, or using the Centroid method. The X-axis represents the number of nodes used to descretize the unit sphere.}
    \label{fig: sphere error}
\end{figure}
By descretizing the unit sphere into triangles, one can integrate $K_{\mathbb{R}^3}(\x, \cdot)$ over the unit sphere by summing the integral over all triangles. 
Since we are evaluating \autoref{eq: integral of K on sphere} on the vertices of each triangle, we encounter a lot of singular integrals. 
As shown in \autoref{fig: sphere error}, the effect of using Quadratic Surface Approximation on these singular integrals significantly reduces the relative error compared to setting these singular integrals to zero or by using a simple centroid quadrature rule. 
It also reduces the number of nodes needed to get below the 1\% relative error threshold. 
Quadratic Surface Approximation only needs around 7200 nodes (14400 triangles) to get less than 1\% relative error, while the other two methods require 28800 nodes (57728 triangles). 
This is important because in BEM, one is often solving a system of integral equations. 
A 7000 $\times$ 7000 system is much easier to solve compared to a 28,000 $\times$ 28,000 system. 
Other shapes such as a torus can also be considered, but there is no true answer to compare to.
Even \texttt{integral2} is unable to integrate $K_{\mathbb{R}^3}(\x, \cdot)$ over a torus.

\subsection{Example uses with BEM}

\begin{figure}
    \centering
    \includegraphics[width=0.9\linewidth]{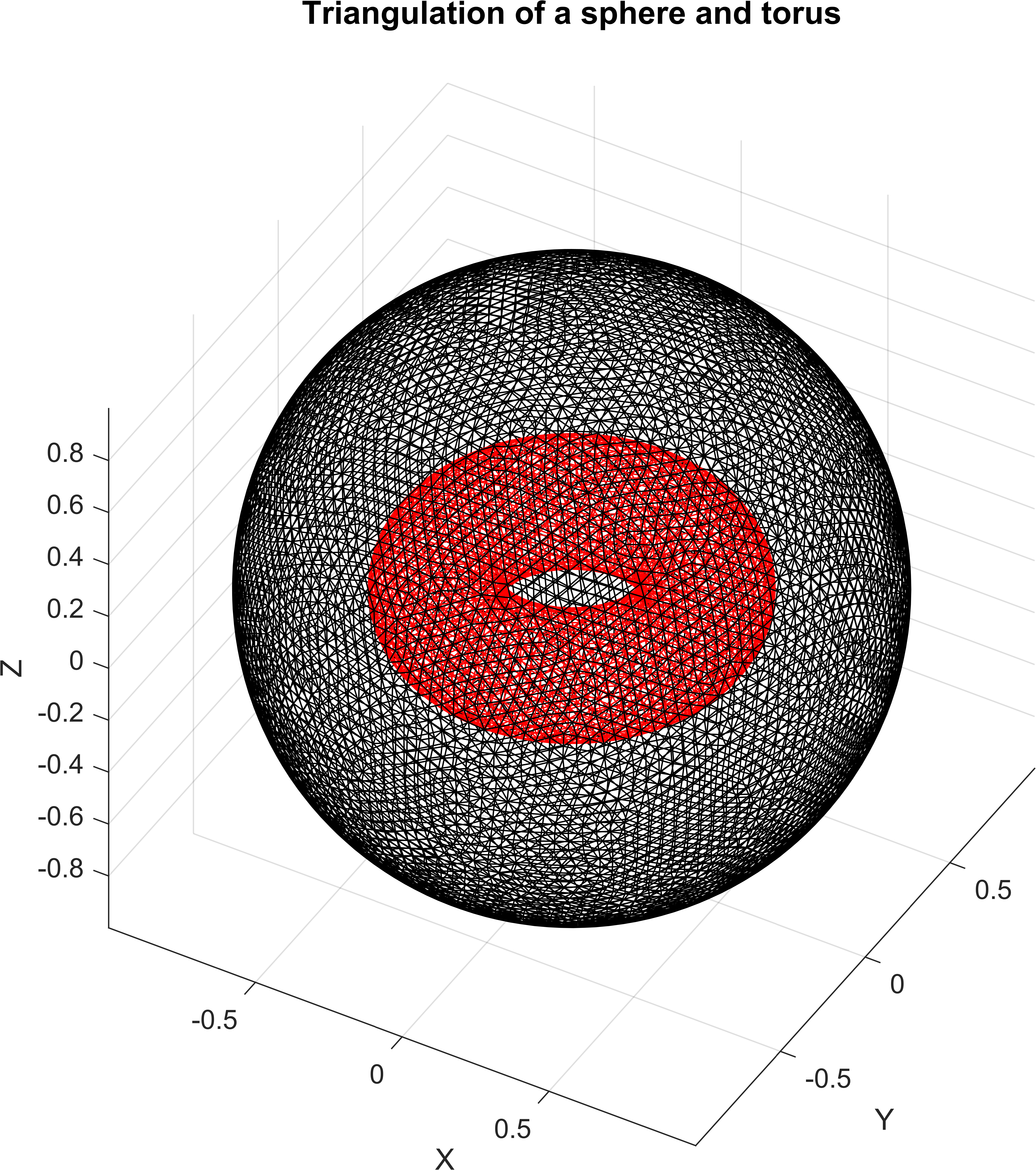}
    \caption{A triangulation of a trous inside of a sphere. Both are centered at the origin. The sphere has radius 1 and the torus has major radius 0.4 and minor radius 0.2. The triangulation of the sphere is shown in black. The triangulation of the torus is shown in red.}
    \label{fig: sphere torus mesh}
\end{figure}

In the next example, let $T$ be a torus aligned with the $XY$-plane with major radius 0.4 and minor radius 0.2.
It is centered at the origin, so it lies within the unit ball $B$. 
Now consider the PDE
\begin{align}
    \nabla \cdot \nabla u(\x) &= 0, \quad \x\in B \backslash \partial T,\\
    \frac{\partial u}{\partial \n} (\x) &= b(\x), \quad \x\in S^2,\\
    \frac{\partial u}{\partial \n^-}(\x) &= 0, x\in \partial T, 
\end{align}
in which $\n^-$ represents the normal that points inwards to the torus.
This represents the potential $u$ inside the ball $B$ with electrical conductivity 1 that also contains a perfect insulating torus. Existence and uniqueness up to a constant of this PDE is discussed in \cite{bower2022fast}.
The Neumann boundary condition $b$ is defined to be 
\begin{equation}
    b(\x) = \cos\left(\arctan\left(\frac{x_2}{x_1}\right)\right).
\end{equation}
In spherical coordinates, the right hand side is simply $\cos(\phi)$, or that current flows in from one side of the ball and leaves on the other side of the ball. 
To use a BEM approach, we first write the solution $u$ as a single layer potential to express the problem as an integral equation.
This integral equation descretized over the flat triangular elements contains integrals of the form \autoref{eq: general integral problem gradient green} with $\x$ being on the triangles. 
The triangulation elements are shown in \autoref{fig: sphere torus mesh}.
Thus, accurately solving this integral equation requires accurate evaluations of singular integrals. 
As the torus and sphere can be defined using signed-distance functions, we can evaluate the normal vectors analytically and use \autoref{eq: second fundamental form} to compute the second fundamental forms.
\begin{figure}
    \centering
    \includegraphics[width=0.9\linewidth]{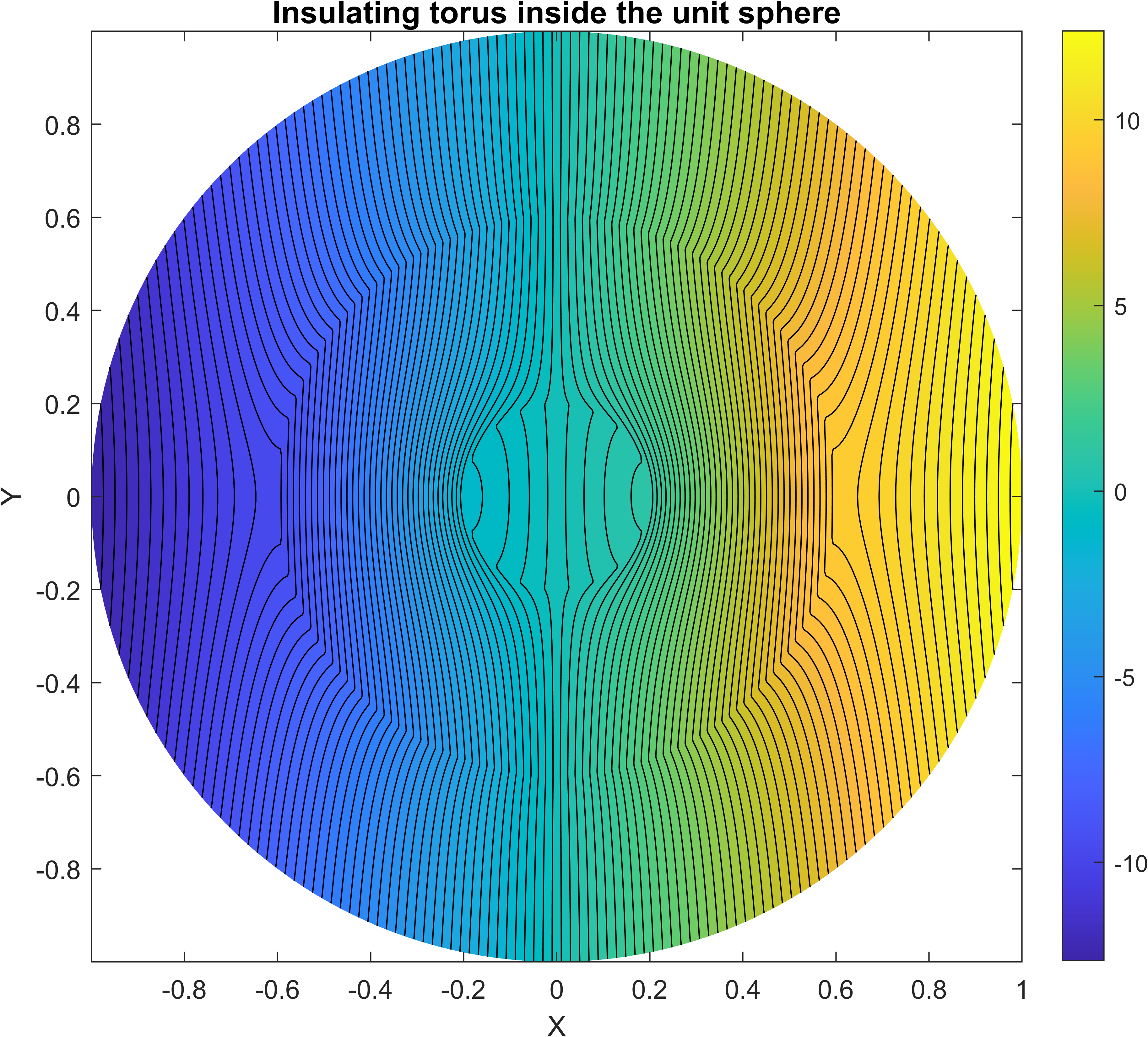}
    \caption{Contour lines of the potential $u$ on the slice $z=0$ inside a conducting ball that has a insulating torus at its center. As the torus is a perfect insulator, the contour lines outside of it is perpendicular to the torus boundary. This allows us to clearly see the location of the torus.}
    \label{fig: circle torus contour}
\end{figure}
As the torus is an insulator, current from outside the torus must flow tangentially to its boundary. 
Since current is the gradient of potential, the contour lines must be perpendicular to the boundary of the torus from the outside \cite{griffith2013electrophysiology}.
As shown in \autoref{fig: circle torus contour}, with accurate integrals, we can very visibly see the location of the insulating torus. 

\begin{figure}
\begin{subfigure}[t]{0.3\textwidth}
    \centering
    \includegraphics[width=\textwidth]{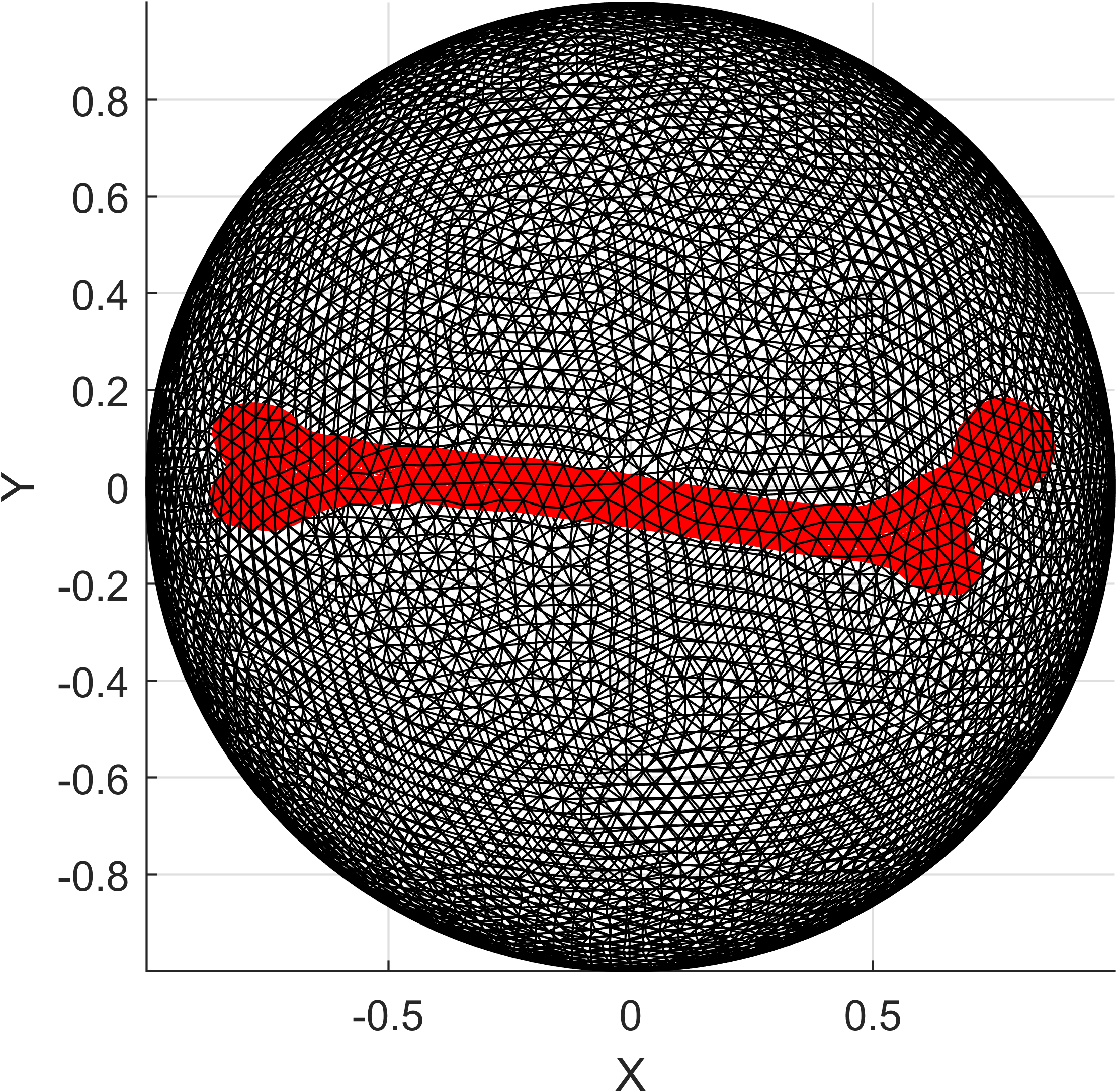}
    \caption{Femur is placed horizontally}
\end{subfigure}
\begin{subfigure}[t]{0.3\textwidth}
    \centering
    \includegraphics[width=\textwidth]{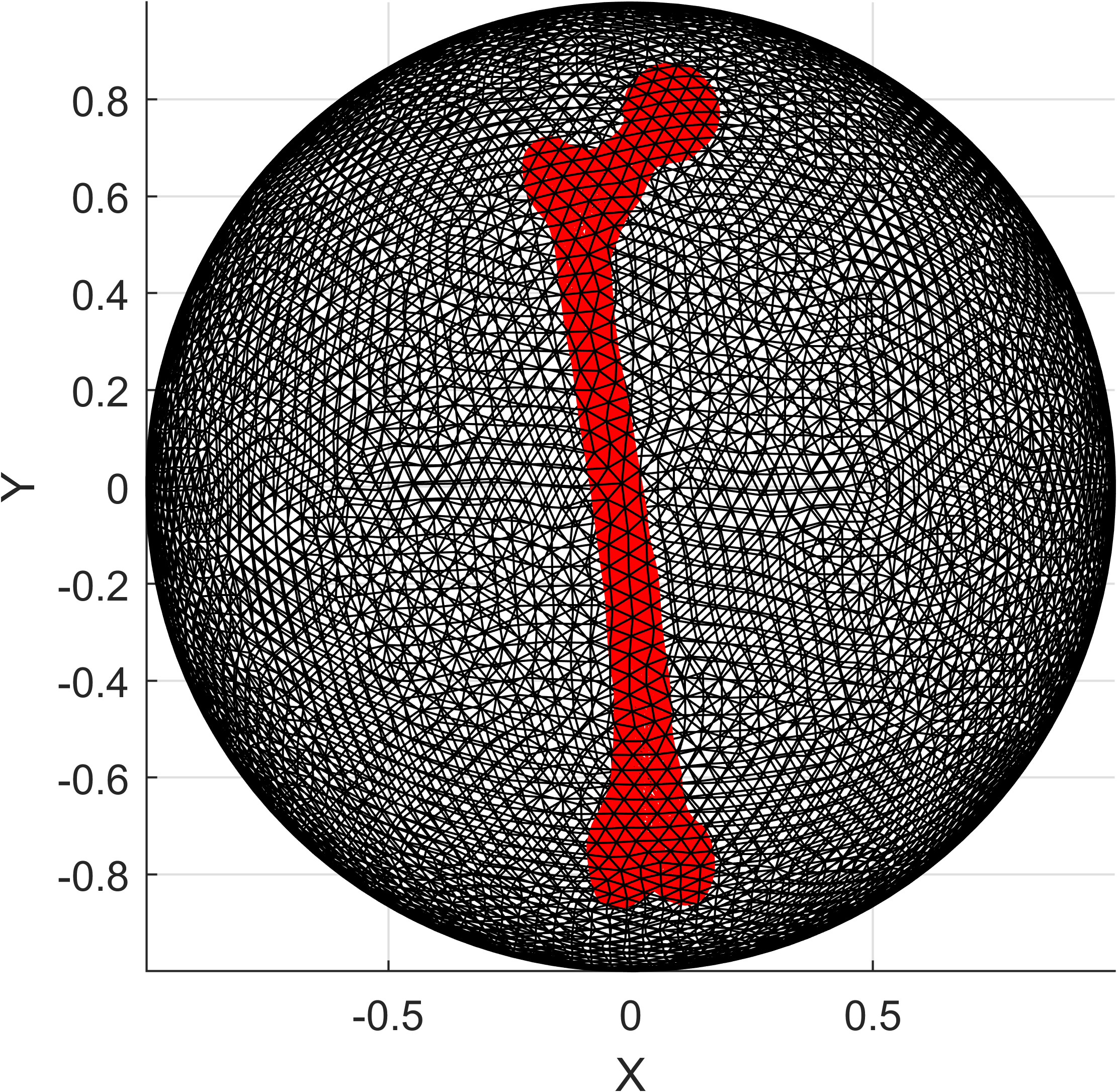}
    \caption{Femur is placed vertically}
\end{subfigure}
\begin{subfigure}[t]{0.3\textwidth}
    \centering
    \includegraphics[width=\textwidth]{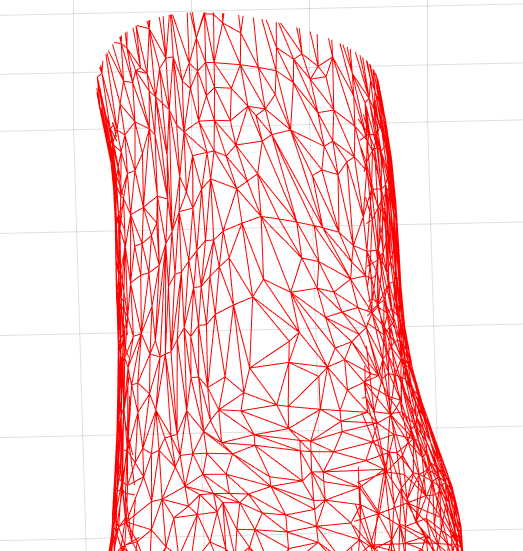}
    \caption{Close up of femur mesh}
\end{subfigure}
    \caption{A triangular mesh of a human femur inside of the unit sphere aligned to the X-axis (a), the Y-axis (b), and the close up of the Femur mesh (c). The sphere mesh is comprised of 7218 nodes (14432 triangles) and the femur mesh is comprised of 8916 nodes (17828 triangles). The triangulation of the sphere is shown in black. The triangulation of the femur is shown in red. The femur mesh has many irregularly shaped triangles.}
    \label{fig: femur mesh}
\end{figure}

\begin{figure}
\begin{subfigure}[t]{0.48\textwidth}
    \centering
    \includegraphics[width=\textwidth]{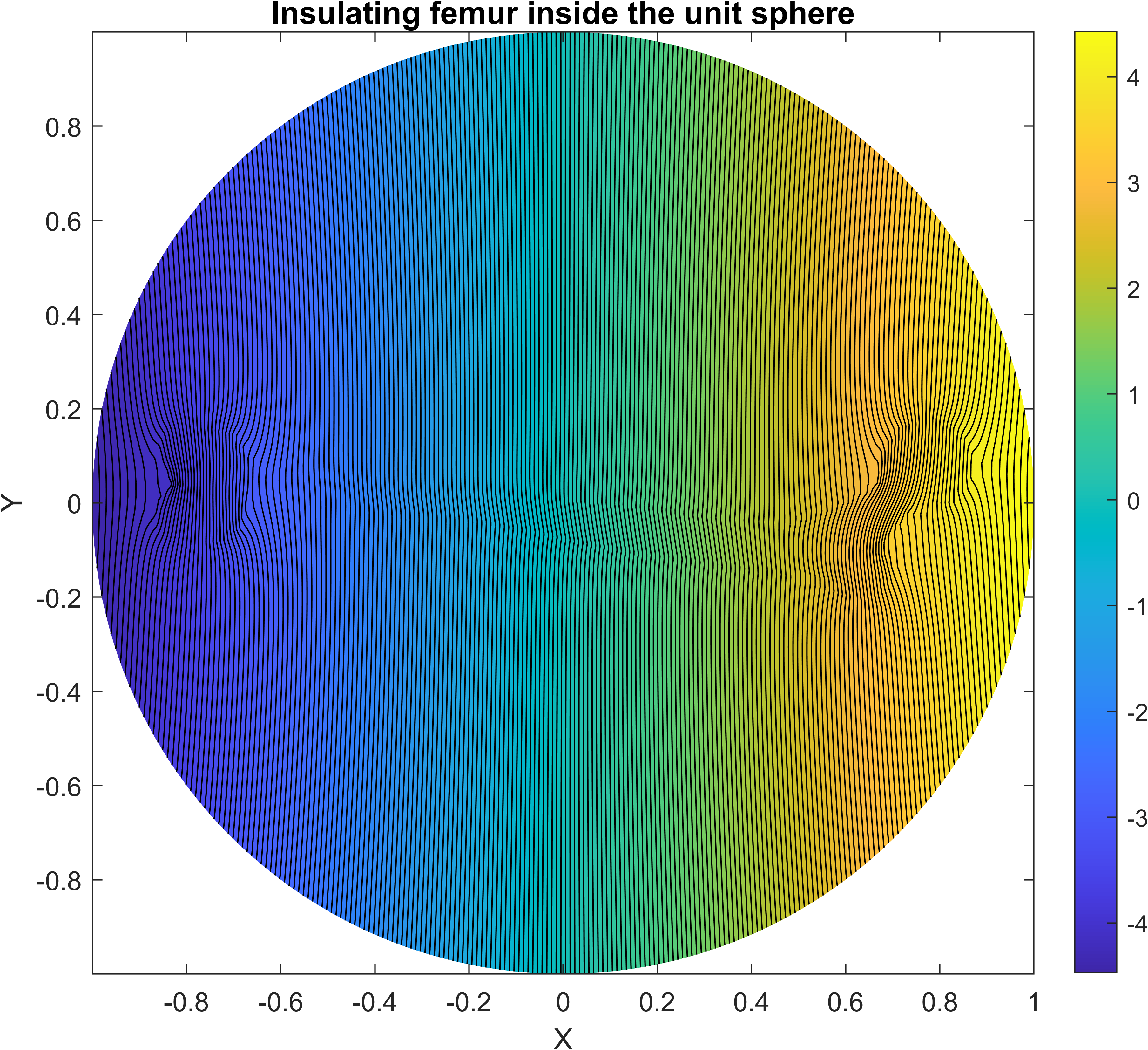}
    \caption{Femur is placed horizontally}
\end{subfigure}
\begin{subfigure}[t]{0.48\textwidth}
    \centering
    \includegraphics[width=\textwidth]{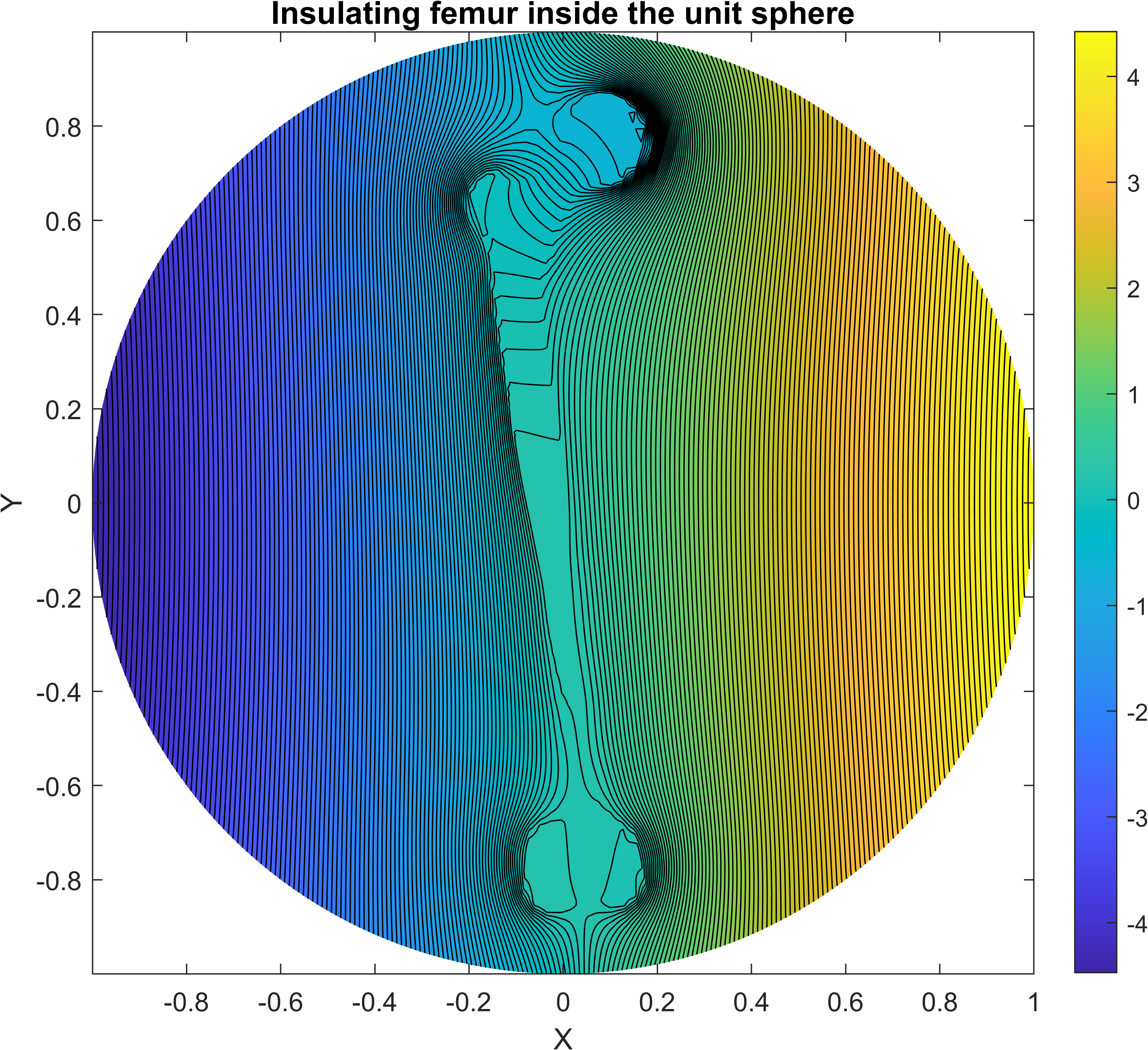}
    \caption{Femur is placed vertically}
\end{subfigure}
\caption{Contour plots of the potential $u$ inside the unit sphere with a femur aligned with the X-axis (a) or the Y-axis (b). In both cases, the contour lines are orthogonal to the boundary of the insulating femur bone.}
\label{fig: contour lines of bone}
\end{figure}

If the user was more interested in human anatomy, they may consider bones instead of a torus. In the next example, we have the same PDE, but we switch the torus with the femur bone of a human.
A triangular mesh of the human femur bone was obtained from the 3D printing website Thingsly \cite{magicbunnydesigns2023} and then coarsened to get a triangular mesh of 8916 nodes.
As shown in \autoref{fig: femur mesh}, the geometry of this new shape is much more complicated and cannot be represented with a signed distance function. Thus, we need to compute the second fundamental form using polynomial approximations as explained in \autoref{sec: approximate second fundamental form using polynomials}. The irregular shaped triangles in the mesh also makes this an ill-conditioned problem to solve with finite element methods or simple quadrature methods. 
Similar to the torus simulation, \autoref{fig: contour lines of bone} show that the contour lines of the potential $u$ are orthogonal to the insulating inner body. 

It is hard to decern accuracy of the solution $u$ as there is no analytic solution to even simple geometries, but the use of Quadratic Surface Approximation and analytic formulas result in more accurate integral evaluations which then allows the linear system to have a smaller condition number.  


\section{Conclusion}\label{sec: conclusion}

In the fields of numerical PDEs, many problems are reformulated into integral equations. For Fredholm integral equations, the kernel functions are often weakly singular on the true geometry, and becomes strongly singular when the domain of integration is approximated by a finite number of elements. This is often the case in BEM, so it is vital to be able to calculate these integrals accurately. Previous practices include setting the strongly singular integrals to zero and stating the error is bounded by the size of the simplex, evaluating only the non-singular part of the integral, or using adaptive refinement with a combination of the previous methods. Most previous analytic methods were only applicable to the 2D case, or when the singularity was of type $r^n, n\geq -1$.

However, with our new method, we can provide fast computational algorithms that lowers the error of these strongly singular integrals even the case in which the kernel function $K(\x, \y)$ has singularity $r^{-2}$.
Quadtratic Approximation incorporates the true geometry of $\partial\Omega$ in a push-forward map so that the singularity caused by approximating $\partial\Omega$ by $\Delta$ is never present in the first place.
This is better than simply setting the integral to zero and much faster than adaptive refinement methods.
However, Quadratic Surface Approximation requires an extra step of calculating the second fundamental form, which can either be calculated with a function parameterization of $\partial\Omega$, or by solving a least squares problem using mesh information. 
Even when the integral is not strongly singular but ``near-singular" due to $\x$ being near $\Delta$, we provide analytic formulas that are just as accurate as adaptive methods when integrating on the simplex $\Delta$ while being much faster. 
The methodology of determining the bounds on integration on an arbitrary flat triangle can also be used to integrate any type of function. 

All formulas for Quadratic Surface Approximation for singular integrals and the analytic method for near singular methods are presented for when $p(\y)$ is a degree two polynomial of its components, but can easily be extended to higher polynomials.
Analytic formulas for integrating $G(\x, \y)p(\y)$ are also presented for linear $p$. 
Since $G$ has a singularity of type $r^{-1}$, it is weakly singular even on $\Delta$, so they are accurate no matter where $\x$ resides.

We believe that with these new formulas, many numerical algorithms of solving PDEs such as BEM \cite{ramvsak20073d, ren2015analytical, bohm2024efficient, atkinson1997numerical} can greatly increase their accuracy without the need to heavily refine their discretization of the domain. As BEM requires solving a dense linear system, not needing to refine triangular elements allows for smaller matrices, faster algorithms, and also smaller memory usage on a computer.
We also note that our methods of integrating singular and near singular integrals can be used in fast multipole methods, as it only affects the near interactions. 


\section*{Acknowledgements}
The authors would like to thank Kirill Serkh for many helpful discussions.

\section*{CRediT authorship contribution statement}

\textbf{Andrew Zheng}: Formal analysis, Software, Writing - original draft.
\textbf{Spyros Alexakis}: Conceptualization, Writing - review and editing. 
\textbf{Adam R. Stinchcombe}: Conceptualization, Supervision, Writing - review and editing.

\section*{Declaration of competing interest}
The authors have no competing interests to declare that are relevant to the content of this article.

\section*{Data availability}
No data was used for the research described in the article. The Rust code implementation is available at \url{https://github.com/AfZheng126/AnalyticIntegrals}.

\section*{Funding}
This work was supported by the Natural Sciences and Engineering Research Council of Canada (NSERC) (Grant numbers RGPIN-2020-01113 and RGPIN-2019-06946).

\appendix
\section{Formulas for near singular integrals}\label{appendix: formula for integral of normal derivative of single layer potential}

The integrals that we want to analytically write out from \autoref{sec: near singular integrals} in polar coordinates are 
\begin{equation}\label{single-layer potential on triangle}
    \int_{r_\mathrm{start}}^{r_\mathrm{end}} \int_{\theta_\mathrm{start}(r)}^{\theta_\mathrm{end}(r)} \frac{r^{a+b+1}\cos^a(\theta)\sin^b(\theta)}{(r^2+c^2)^{3/2}} ~\diff \theta \diff r,
\end{equation}
in which $(a,b)\in \left\{(0,0),(1,0),(0,1),(2,0),(0,2),(1,1)\right\}$. The integrals in $\theta$ are easy to integrate. After integrating, we substitute $\theta = \phi \pm \arccos(d/r)$ to see what we need to integrate in $r$. This series of steps gives us

{\small
\begin{align*}
(a, b) = (0, 0): &\quad \phi \pm \arccos(d/r), \\
(a, b) = (1, 0): &\quad \frac{d}{r} \sin\phi \pm \frac{\sqrt{r^2-d^2}}{r}\cos\phi, \\
(a, b) = (0, 1): &\quad -\frac{d}{r} \cos\phi \pm \frac{\sqrt{r^2-d^2}}{r}\sin\phi, \\
(a, b) = (2, 0): &\quad \frac12 (\phi \pm \arccos(d/r)) \pm \cos(2\phi)\frac{d\sqrt{r^2-d^2}}{2r^2} + \left(\frac{d^2}{2r^2}-\frac14\right)\sin(2\phi), \\
(a, b) = (0, 2): &\quad \frac12 (\phi \pm \arccos(d/r)) \mp \cos(2\phi)\frac{d\sqrt{r^2-d^2}}{2r^2} - \left(\frac{d^2}{2r^2}-\frac14\right)\sin(2\phi), \\
(a, b) = (1, 1): &\quad -\frac{d^2}{2r^2}\cos(2\phi) - \frac{1}{2}(\sin(\phi))^2 \pm \frac{d\sqrt{r^2-d^2}}{2r^2}\sin(2\phi).\\
\end{align*}
}%
After substitution of these $\theta$ integrals into \autoref{single-layer potential on triangle}, there are five integrals in $r$,
{\small
\begin{gather*}
    \int \frac{r}{(r^2+c^2)^{3/2}}~ \diff r = -\frac{1}{\sqrt{r^2+c^2}},\\
    \int \frac{r}{(r^2+c^2)^{3/2}} \arccos{\frac{d}{r}}~\diff r = \begin{cases}
        -\frac{\arccos{\frac{d}{r}}}{\sqrt{r^2+c^2}}+\frac{1}{c}\arctan\left(\frac{c}{d}\sqrt{\frac{r^2-d^2}{r^2+c^2}}\right), & d \ne 0, r\neq 0 \\
        -\frac{\pi/2}{\sqrt{r^2+c^2}}, & d = 0 
    \end{cases},\\
    \int \frac{r}{(r^2+c^2)^{3/2}}\sqrt{r^2-d^2}~\diff r = \frac12 \log \left(2 \sqrt{r^2+c^2}
    \sqrt{r^2-d^2}+c^2-d^2+2 r^2\right)-\sqrt{\frac{r^2-d^2}{r^2+c^2}},\\
    \int \frac{r^3}{(r^2+c^2)^{3/2}}~\diff r = \frac{r^2 + 2c^2}{\sqrt{r^2+c^2}},\\
\end{gather*}
\vspace{-1cm} 
\begin{align*}
    &\int \frac{r^3}{(r^2+c^2)^{3/2}} \arccos{\frac{d}{r}}~\diff r =\\
        &\qquad \begin{cases}
            -2 c~\text{arctan} \left(\frac{c}{d} \sqrt{\frac{r^2-d^2}{r^2+c^2}} \right)+ \frac{r^2 + 2c^2}{\sqrt{r^2+c^2}}\arccos\frac{d}{r} - d~ \text{arctanh}\left( \sqrt{\frac{r^2-d^2}{r^2+c^2}} \right), & d \ne 0, r\ne 0 \\
            \frac{\pi}{2}\frac{r^2 + 2c^2}{\sqrt{r^2+c^2}}, & d = 0
        \end{cases}.
\end{align*}
}%
We denote these integrals as $R_1, R_2, \ldots, R_5$ respectively.
Due to the first integral, we see that $c$ cannot equal to 0 when $r=0$ (see \autoref{equation: I1}).
In the case when $c=0$, the five integrals become 
{\small
\begin{gather*}
    R_1 = \int \frac{1}{r^2}~\diff r = -\frac{1}{r},\\
    R_2 = \int \frac{1}{r^2} \arccos{\frac{d}{r}}~\diff r = \begin{cases}
        -\frac{\arccos\frac{d}{r}}{r} + \frac{\sqrt{r^2-d^2}}{rd}, & d \ne 0, l\neq 0 \\
        -\frac{\pi/2}{r}, & d = 0 
    \end{cases},\\
    R_3 = \int \frac{1}{r^2}\sqrt{r^2-d^2}~\diff r = 
        \begin{cases}
            -\frac{\sqrt{r^2 - d^2}}{r} + \log \left( \frac{\sqrt{r^2-d^2} + r}{d}\right), & d\ne 0 \\
            \log(r), & d = 0
        \end{cases},\\
    R_4 = \int 1~\diff r = r, \\
    R_5 = \int \arccos{\frac{d}{r}}~\diff r =
    \begin{cases}
        r \arccos(\frac{d}{r}) - d \log\left( \frac{\sqrt{r^2 - d^2} + r}{d}\right), & d \ne 0, r\ne 0 \\
        \frac{\pi}{2}r, & d = 0
    \end{cases}.
\end{gather*}
}%
Due to numerical precision, some terms can give NaN values (such as $\text{arctanh}\left( \sqrt{\frac{r^2-d^2}{r^2+c^2}} \right)$ when $d^2, c^2 \approx 10^{-16}$), so these terms are calculated using the first four terms of their Taylor expansion. 
Now we can write out the full formulas for the integrals. For the sake of notation, we write 
\begin{equation*}
    R_2(a, b, c, d) = \int_a^b \frac{r}{(r^2+c^2)^{3/2}} \arccos{\frac{d}{r}}~\diff r
\end{equation*}
and similar for other integrals $R_i$ except for $R_1$ and $R_4$ which does not depend on $d$.

For $I_0$, we have 
\begin{equation}\label{equation: I1}
\begin{split}
    I_0  & = \int_{r_\mathrm{start}}^{r_\mathrm{end}} \int_{\phi_\mathrm{start} + \mathrm{sign}_\mathrm{start}\arccos\left(d_\mathrm{start}/r\right)}^{\phi_\mathrm{end} + \mathrm{sign}_\mathrm{end}\arccos\left(d_\mathrm{end}/r\right)} \frac{r}{(r^2 + c^2)^{3/2}}\, \diff \theta \, \diff r \\
    & = (\phi_\mathrm{end} - \phi_\mathrm{start}) R_1(r_\mathrm{start}, r_\mathrm{end}, c) \\
    & + \mathrm{sign}_\mathrm{end} R_2(r_\mathrm{start}, r_\mathrm{end}, c, d_\mathrm{end}) \\
    & - \mathrm{sign}_\mathrm{start} R_2(r_\mathrm{start}, r_\mathrm{end}, c, d_\mathrm{start}).
\end{split}
\end{equation}
It is important to know that when we have $\phi_\mathrm{end} - \phi_\mathrm{start}$, we need to check for branch cuts in $\arctan$. 
These is because $\phi$ was calculated via $\arctan$, which is multi-valued. 
When checking the branch cut, we just need to check if $\theta_\mathrm{end} - \theta_\mathrm{start}\in [0, 2\pi]$, not if $\phi_\mathrm{end} - \phi_\mathrm{start}\in [0, 2\pi]$. 
This branch cut problem also shows up in $I_{x^2}$ and $I_{y^2}$. 
We now list the other integrals.

{\small
\begin{equation}
    \begin{split}
        I_x & = \left(d_\mathrm{end}\sin(\phi_\mathrm{end}) - d_\mathrm{start}\sin(\phi_\mathrm{start})\right) R_1(r_\mathrm{start}, r_\mathrm{end}, c) \\
        & + \mathrm{sign}_\mathrm{end}\cos(\phi_\mathrm{end}) R_3(r_\mathrm{start}, r_\mathrm{end}, c, d_\mathrm{end}) \\
        & - \mathrm{sign}_\mathrm{start}\cos(\phi_\mathrm{start}) R_3(r_\mathrm{start}, r_\mathrm{end}, c, d_\mathrm{start}).
    \end{split}
\end{equation}
\begin{equation}
\begin{split}
    I_y & = \left(d_\mathrm{start}\cos(\phi_\mathrm{start}) - d_\mathrm{end}\cos(\phi_\mathrm{end})\right)  R_1(r_\mathrm{start}, r_\mathrm{end}, c)\\
    & + \mathrm{sign}_\mathrm{end}\sin(\phi_\mathrm{end}) R_3(r_\mathrm{start}, r_\mathrm{end}, c, d_\mathrm{end}) \\
    & - \mathrm{sign}_\mathrm{start}\sin(\phi_\mathrm{start}) R_3(r_\mathrm{start}, r_\mathrm{end}, c, d_\mathrm{start}).
\end{split}
\end{equation}
\begin{equation}
    \begin{split}
        I_{x^2} & = \frac{1}{2}\left(\phi_\mathrm{end} - \frac{\sin(2\phi_\mathrm{end})}{2} - \phi_\mathrm{start} + \frac{\sin(2\phi_\mathrm{start})}{2}\right)R_4(r_\mathrm{start}, r_\mathrm{end}, c) \\
        & + \frac{\mathrm{sign}_\mathrm{end}}{2}R_5(r_\mathrm{start}, r_\mathrm{end}, c, d_\mathrm{end})\\
        & - \frac{\mathrm{sign}_\mathrm{start}}{2}R_5(r_\mathrm{start}, r_\mathrm{end}, c, d_\mathrm{start})\\
        & + \frac{1}{2}(d_\mathrm{end}^2\sin(2\phi_\mathrm{end}) - d_\mathrm{start}^2\sin(2\phi_\mathrm{start})) R_1(r_\mathrm{start}, r_\mathrm{end}, c)\\
        & + \frac{1}{2}\mathrm{sign}_\mathrm{end}d_\mathrm{end}\cos(2\phi_\mathrm{end}) R_3(r_\mathrm{start}, r_\mathrm{end}, c, d_\mathrm{end})\\
        & - \frac{1}{2}\mathrm{sign}_\mathrm{start}d_\mathrm{start}\cos(2\phi_\mathrm{start}) R_3(r_\mathrm{start}, r_\mathrm{end}, c, d_\mathrm{start}).
    \end{split}
\end{equation}
\begin{equation}
    \begin{split}
        I_{y^2} & = \frac{1}{2}\left(\phi_\mathrm{end} + \frac{\sin(2\phi_\mathrm{end})}{2} - \phi_\mathrm{start} - \frac{\sin(2\phi_\mathrm{start})}{2}\right)R_4(r_\mathrm{start}, r_\mathrm{end}, c) \\
        & + \frac{\mathrm{sign}_\mathrm{end}}{2}R_5(r_\mathrm{start}, r_\mathrm{end}, c, d_\mathrm{end})\\
        & - \frac{\mathrm{sign}_\mathrm{start}}{2}R_5(r_\mathrm{start}, r_\mathrm{end}, c, d_\mathrm{start})\\
        & + \frac{1}{2}(- d_\mathrm{end}^2\sin(2\phi_\mathrm{end}) + d_\mathrm{start}^2\sin(2\phi_\mathrm{start})) R_1(r_\mathrm{start}, r_\mathrm{end}, c)\\
        & - \frac{1}{2}\mathrm{sign}_\mathrm{end}d_\mathrm{end}\cos(2\phi_\mathrm{end}) R_3(r_\mathrm{start}, r_\mathrm{end}, c, d_\mathrm{end})\\
        & + \frac{1}{2}\mathrm{sign}_\mathrm{start}d_\mathrm{start}\cos(2\phi_\mathrm{start}) R_3(r_\mathrm{start}, r_\mathrm{end}, c, d_\mathrm{start}).
    \end{split}
\end{equation}
\begin{equation}
    \begin{split}
        I_{xy} & = \frac{1}{2}(\sin^2(\phi_\mathrm{start}) - \sin^2(\phi_\mathrm{end})) R_4(r_\mathrm{start}, r_\mathrm{end}, c)\\
        & + \frac{1}{2}(d_\mathrm{start}^2 \cos(2\phi_\mathrm{start}) - d_\mathrm{end}^2\cos(2\phi_\mathrm{end})) R_1(r_\mathrm{start}, r_\mathrm{end}, c)\\
        & + \frac{1}{2}\mathrm{sign}_\mathrm{end}\sin(2\phi_\mathrm{end})d_\mathrm{end} R_3(r_\mathrm{start}, r_\mathrm{end}, c, d_\mathrm{end})\\
        & - \frac{1}{2}\mathrm{sign}_\mathrm{start}\sin(2\phi_\mathrm{start})d_\mathrm{start} R_3(r_\mathrm{start}, r_\mathrm{end}, c, d_\mathrm{start}).
    \end{split}
\end{equation}
}%
For higher powers, integrating in theta first gives
{\small
\begin{align*}
(a, b) = (3, 0): &\quad \int \cos^3(\theta)\, \diff \theta = \sin(\theta) - \frac{1}{3}\sin^3(\theta),\\
(a, b) = (2, 1): &\quad \int \cos^2(\theta)\sin(\theta)\, \diff \theta = -\frac{1}{3} \cos^3(\theta),\\
(a, b) = (1, 2): &\quad \int \cos(\theta)\sin^2(\theta)\, \diff \theta = \frac{1}{3} \sin^3(\theta),\\
(a, b) = (0, 3): &\quad \int \sin^3(\theta)\, \diff \theta = -\cos(\theta) + \frac{1}{3}\cos^3(\theta).
\end{align*}
}%
Plugging in $\theta = \phi\pm\arccos(d/r)$, we get for $(a, b) = (3, 0)$
\begin{equation*}
\begin{split}
    (a, b) = (3, 0) &:\left(\sin(\phi)\frac{d}{r} \pm \cos(\phi) \frac{\sqrt{r^2 - d^2}}{r}\right) \\
&-\frac{1}{3}\left[ \sin^3(\phi)\frac{d^3}{r^3} \pm 3 \sin^2(\phi)\cos(\phi) \frac{d^2\sqrt{r^2-d^2}}{r^3}\right. \\ 
&+ \left. 3\sin(\phi)\cos^2(\phi)\frac{d(r^2-d^2)}{r^3} \pm \cos^3(\phi) \frac{(r^2-d^2)^{3/2}}{r^3}\right].
\end{split}    
\end{equation*}
For $(a, b) = (0, 3)$, we have 
\begin{equation*}
\begin{split}
    (a, b) = (0, 3) &:-\left(\cos(\phi)\frac{d}{r} \mp \sin(\phi) \frac{\sqrt{r^2 - d^2}}{r}\right) \\
&+\frac{1}{3}\left[ \cos^3(\phi)\frac{d^3}{r^3} \mp 3 \cos^2(\phi)\sin(\phi) \frac{d^2\sqrt{r^2-d^2}}{r^3}\right. \\ 
&+ \left. 3\cos(\phi)\sin^2(\phi)\frac{d(r^2-d^2)}{r^3} \mp \sin^3(\phi) \frac{(r^2-d^2)^{3/2}}{r^3}\right].
\end{split}    
\end{equation*}
Integrating in $r$, the new integral we have to evaluate is
\begin{equation*}
    \begin{split}
        \int \frac{r(r^2-d^2)^{3/2}}{(r^2+c^2)^{3/2}}\, \diff r &= 
        \sqrt{\frac{r^2-d^2}{r^2+c^2}} \frac{r^2+2d^2+3c^2}{2} \\
        &- \frac{3(c^2+d^2)\log(\sqrt{r^2+c^2} + \sqrt{r^2-d^2})}{2}.
    \end{split}
\end{equation*}
When $c=0$, we have that the integral becomes
\begin{equation*}
    \int \frac{(r^2-d^2)^{3/2}}{r^2}\, \diff r = \left(\frac{d^2}{r} + \frac{r}{2}\right)\sqrt{r^2-d^2} - \frac{3d^2}{2}\log(r + \sqrt{r^2-d^2}).
\end{equation*}
These equations can then be used to get the integrals for $(a, b) = (2, 1)$ or $(a, b) = (1, 2)$.

\section{Formulas for Quadratic Surface Approximation when the target is a vertex}\label{appendix: formula for quadratic approximation}
Recall that in Quadratic Surface Approximation when $\x$ is a vertex of the triangle, the integrals are of the form
\begin{equation}
    \int_0^{\theta_\mathrm{end}} (\cos(\theta))^a (\sin(\theta))^b \int_0^{r(\theta)} r^{a+b+1-3}\, \diff r \, \diff \theta.
\end{equation}
This simple formula is only because we can fix one edge to be the $X$-axis, so that the integral in $\theta$ starts at zero. The integrals in $r$ are trivial as $a+b+1-3\geq 0$, so we get 
\begin{equation}
    \int_0^{\theta_\mathrm{end}} (\cos(\theta))^a (\sin(\theta))^b \frac{1}{a+b-1}\left(\frac{|p\V_2|\sin(\theta_2)}{\sin(\theta+\theta_2)}\right)^{a+b-1}\, \diff \theta.
\end{equation}
We shall only consider the cases in which $(a, b)\in \{(1, 1), (2, 0), (0, 2), (2, 1), (1, 2), (3, 0), (0, 3)\}$. These correspond to when we approximate $p$ using linear functions, but it is not difficult to write similar formulas for higher powers. 

Ignoring the constants in the integrals, we present the formulas for integrals of the form 
\begin{equation}
    \int_0^{\theta_\mathrm{end}} (\cos(\theta))^a (\sin(\theta))^b \left(\frac{1}{\sin(\theta+\theta_2)}\right)^{a+b-1}\, \diff \theta.
\end{equation}
{\small
\begin{align*}
    (a, b) = (1, 1) &= \sin(\theta-\theta_2) - \frac{\sin(2\theta_2)}{2} \log\left(\tan\left(\frac{\theta+\theta_2}{2}\right)\right),\\
    (a, b) = (2, 0) &= \cos(\theta - \theta_2) + \cos^2(\theta_2)\log\left(\tan\left(\frac{\theta + \theta_2}{2}\right)\right),\\
    (a, b) = (0, 2) &=  -\cos(\theta - \theta_2) + \sin^2(\theta_2)\log\left(\tan\left(\frac{\theta + \theta_2}{2}\right)\right),\\
    (a, b) = (2, 1) &= \frac{1}{4}\big(-2(\cos(\theta_2) + 3\cos(3\theta_2))\mathrm{arctanh}\left(\cos(\theta_2) - \sin(\theta_2)\tan\left(\frac{\theta}{2}\right)\right) \\
    &\qquad + \csc(\theta + \theta_2)(2\sin(2\theta-\theta_2) + \sin(\theta_2) + 3\sin(3\theta_2))\big),\\
    (a, b) = (1, 2) &= \frac{1}{4}\big(-2(\sin(\theta_2) - 3\sin(3\theta_2))\mathrm{arctanh}\left(\cos(\theta_2) - \sin(\theta_2)\tan\left(\frac{\theta}{2}\right)\right) \\
    &\qquad - \csc(\theta + \theta_2)(2\cos(2\theta-\theta_2) + \cos(\theta_2) - 3\cos(3\theta_2))\big),\\
    (a, b) = (3, 0) &= -\cos(\theta-2\theta_2) - 6\cos(\theta_2)\sin^2(\theta_2)\mathrm{arctanh}\left( \cos(\theta_2) - \sin(\theta_2)\tan\left(\frac{\theta}{2}\right)\right)\\
    &\qquad + \sin^3(\theta_2)\csc(\theta+\theta_2),\\
    (a, b) = (0, 3) &= -\sin(\theta-2\theta_2) - 6\cos^2(\theta_2)\sin(\theta_2)\mathrm{arctanh}\left(\cos(\theta_2) - \sin(\theta_2)\tan\left(\frac{\theta}{2}\right)\right)\\
    &\qquad -\cos^3(\theta_2)\csc(\theta+\theta_2).
\end{align*}
}%
\section{Formulas for the Green's function}\label{appendix: formula for single layer potential}
Though most of the paper discusses integrating $K(\x, \y) p(\y)$, we also provide analytic equations for integrating $G(\x, \y)p(\y)$ where $G$ is defined from \autoref{eq: Green in 3D} and $p(\y)$ is a linear polynomial.
Given 
{\small
\begin{gather*}
J_0:=\int_\Delta \frac{1}{(y_x^2+y_y^2+c^2)^{1/2}} ~\diff A_{\y}, \quad
J_x:=\int_\Delta \frac{y_x}{(y_x^2+y_y^2+c^2)^{1/2}}~\diff A_{\y}, \\
J_y:=\int_\Delta \frac{y_y}{(y_x^2+y_y^2+c^2)^{1/2}}~\diff A_{\y}.
\end{gather*}
}%
The desired integral is $J = -\frac{1}{4\pi}\sum_{i=1}^3 \gamma(\V_i) \big[ l_{i,0}J_0 + l_{i,x}J_x + l_{i,y}J_y\big].$
The integrals are of the form
\begin{equation}
    \int_{r_\mathrm{start}}^{r_l} \int_{\theta_\mathrm{start}}^{\theta_\mathrm{end}} \frac{r^{a+b+1}(\cos(\theta))^a(\sin(\theta))^b}{(r^2 + c^2)^{1/2}}~\diff \theta\diff r
\end{equation}
for $(a, b) = \{ (0, 0), (0, 1), (1, 0)\}$. The integrals in $\theta$ are 
{\small
\begin{align*}
    (a, b) = (0, 0): &\quad \int 1 \, \diff \theta = \theta = \phi \pm \arccos\left(\frac{d}{r}\right), \\
    (a, b) = (1, 0): &\quad \int \cos(\theta) \, \diff \theta = \sin(\theta) = \sin\left(\phi \pm \arccos\left(\frac{d}{r}\right)\right) \\
    &\qquad =  \frac{d}{r} \sin(\phi) \pm \frac{\sqrt{r^2-d^2}}{r}\cos(\phi),\\
    (a, b) = (0, 1): &\quad \int \sin(\theta)\, \diff \theta = -\cos(\theta) = -\cos\left(\phi \pm \arccos\left(\frac{d}{r}\right)\right) \\
    &\qquad =  -\frac{d}{r} \cos(\phi) \pm \frac{\sqrt{r^2-d^2}}{r}\sin(\phi).
\end{align*}
}%
The integrals in $r$, which we denote as $J_1, J_2, J_3$, are respectively 
{\small
\begin{gather*}
    \int \frac{r}{(r^2 + c^2)^{1/2}}\, \diff r = \sqrt{r^2 + c^2}, \\
    \int \frac{r \arccos(\frac{d}{r})}{(r^2 + c^2)^{1/2}}\, \diff r = \begin{cases}
            \sqrt{r^2 + c^2}\arccos(\frac{d}{r}) - c\arcsin\left(\frac{c\sqrt{\frac{r^2-d^2}{r^2}}}{\sqrt{c^2 + d^2}}\right) - d~\text{arctanh}\left(\sqrt{\frac{r^2-d^2}{r^2+c^2}}\right) & d, r \ne 0\\
            \frac{\pi}{2}\sqrt{r^2 + c^2} & d = 0\\
        \end{cases},\\
    \int \frac{r\sqrt{r^2-d^2}}{(r^2 + c^2)^{1/2}}\, \diff r = \frac{1}{2} \left(\sqrt{r^2+c^2}\sqrt{r^2-d^2} + (c^2+d^2)\log\left(\sqrt{r^2+c^2} - \sqrt{r^2-d^2}\right)\right).
\end{gather*}
}%
The sign of $c$ does not matter on both sides of the equations. In the right hand side of the second integral, $\arcsin$ is an odd function, so the sign of the $c$ inside and in front of $\arcsin$ cancel each other out. $c$ can also equal to 0 on the left hand side even if $r = 0$, which gives us
{\small
\begin{gather*}
    J_1 = \int \frac{r}{(r^2)^{1/2}}\, \diff r = r,\\
    J_2 = \int \frac{r \arccos\left(\frac{d}{r}\right)}{(r^2)^{1/2}}\, \diff r = \begin{cases}
            r\arccos\left(\frac{d}{r}\right) + \frac{d}{2}\log\left(-\frac{\sqrt{r^2-d^2}-r}{\sqrt{r^2-d^2}+r}\right) & d, r \ne 0\\
            \frac{\pi}{2}r & d = 0\\
        \end{cases},\\
    J_3 = \int \frac{r\sqrt{r^2-d^2}}{(r^2)^{1/2}}\, \diff r = \begin{cases}
            \frac{1}{2} \left(r\sqrt{r^2-d^2} + d^2\log\left(r - \sqrt{r^2-d^2}\right)\right) & d, r \ne 0\\
            \frac{1}{2}r^2 & d = 0
        \end{cases}.
\end{gather*}
}%
Hence, with similar notation as before, 
{\small
\begin{align*}
    J_0 & = (\phi_\mathrm{end} - \phi_\mathrm{start}) J_1(r_\mathrm{start}, r_\mathrm{end}, c)\\ 
    & \quad + \mathrm{sign}_\mathrm{end}J_2(r_\mathrm{start}, r_\mathrm{end}, c, d_\mathrm{end}) \\
    & \quad - \mathrm{sign}_\mathrm{start}J_2(r_\mathrm{start}, r_\mathrm{end}, c, d_\mathrm{start}),\\
    J_x  & = (d_\mathrm{end}\sin(\phi_\mathrm{end}) - d_\mathrm{start}\sin(\phi_\mathrm{start})) J_1(r_\mathrm{start}, r_\mathrm{end}, c)\\
    & \quad + \mathrm{sign}_\mathrm{end}\cos(\phi_\mathrm{end})J_3(r_\mathrm{start}, r_\mathrm{end}, c, d_\mathrm{end})\\
    & \quad - \mathrm{sign}_\mathrm{start}\cos(\phi_\mathrm{start})J_3(r_\mathrm{start}, r_\mathrm{end}, c, d_\mathrm{start}),\\
    J_y  & = -(d_\mathrm{end}\cos(\phi_\mathrm{end}) - d_\mathrm{start}\cos(\phi_\mathrm{start})) J_1(r_\mathrm{start}, r_\mathrm{end}, c)\\
    & \quad + \mathrm{sign}_\mathrm{end}\sin(\phi_\mathrm{end})J_3(r_\mathrm{start}, r_\mathrm{end}, c, d_\mathrm{end})\\
    & \quad - \mathrm{sign}_\mathrm{start}\sin(\phi_\mathrm{start})J_3(r_\mathrm{start}, r_\mathrm{end}, c, d_\mathrm{start}).
\end{align*}
}%
Just like before, we need to be careful of the branch cut problem in $J_0$. As all integrals are finite even when $c=0$, $\x$ can be any arbitrary point in $\mathbb{R}^3$.

\section{Integrals of Green and its normal derivative on the unit sphere}\label{appendix: proof of integrals on the unit sphere}

Due to spherical symmetry, without loss of generality let $\x = [0, 0, 1]^T$. 
Then integral can then be explicitly computed using polar coordinates.
\begin{align*}
    \int_{S^2} K(\x, \y)\, dS(\y) &= -\frac{1}{4\pi} \int_0^{2\pi} \int_{0}^{\pi} \frac{(1 - \cos(\theta))\sin(\theta)}{(\sin^2(\theta)\cos^2(\phi) + \sin^2(\theta)\sin^2(\phi) + (1 - \cos(\theta))^2)^{3/2}}\, \diff\theta\diff\phi\\
    & = -\frac{1}{4\pi} \int_0^{2\pi} \int_{0}^{\pi} \frac{(1 - \cos(\theta))\sin(\theta)}{(2(1 - \cos(\theta)))^{3/2}}\, \diff\theta\diff\phi\\
    & = -\frac{1}{2^{1+3/2}} \int_{0}^{\pi} \frac{\sin(\theta)}{\sqrt{1 - \cos(\theta)}}\, \diff\theta\\
    & = -\frac{2\sqrt{2}}{2^{(1+3/2)}}\\
    & = -\frac{1}{2}.
\end{align*}

\bibliographystyle{elsarticle-num} 
\bibliography{bibliography}



\end{document}